\title{Numerical reconstruction of the covariance matrix \\ of a spherically
truncated multinormal distribution\\[2.0ex]}  

\author{
  Filippo Palombi$^{\it a,b}$\footnote{Corresponding author. e--mail: 
  {\tt filippo.palombi$@$enea.it}},\ \ %
  Simona Toti$^{\it a}$ and Romina Filippini$^{\it a}$ \\[4.0ex]
  {$^{\it a}$\it \small ISTAT -- Istituto Nazionale di Statistica}\\[0.5ex]
  {\small Via Cesare Balbo 16, 00184 Rome -- Italy}\\[3.0ex]
  {$^{\it b}$\it \small ENEA -- Italian Agency for New Technologies, Energy}\\[0.0ex]
  {\it \small and Sustainable Economic Development}\\[0.5ex]
  {\small Via Enrico Fermi 45, 00040 Frascati -- Italy}\\[1.0ex]
}
\date{April 2014}

\documentclass[11pt,fleqn]{article}

\usepackage{footmisc,booktabs}
\usepackage{psfrag}
\usepackage{amsmath}
\usepackage{mathtools}
\usepackage{amsfonts}
\usepackage{dsfont}
\usepackage{amsthm}
\usepackage{epsfig}
\usepackage{multicol}
\usepackage[title]{appendix}
\usepackage{pgf}
\usepackage{bbm}
\usepackage[utf8]{inputenc}
\usepackage[margin=30pt,font=small,labelfont=bf,labelsep=endash]{caption}
\usepackage{marvosym}

\numberwithin{equation}{section}\usepackage{amssymb}
\definecolor{Blue}{rgb}{0,0,1}
\definecolor{Black}{rgb}{0,0,0}

\newcommand{\cB}{{\cal B}}
\newcommand{\cC}{{\cal C}}
\newcommand{\cD}{{\cal D}}
\newcommand{\cE}{{\cal E}}

\newcommand{\cN}{{\cal N}}

\newcommand{\cR}{{\cal R}}
\newcommand{\cS}{{\cal S}}
\newcommand{\cW}{{\cal W}}

\newcommand{\fS}{{\frak S}}
\newcommand{\ilr}{{\rm ilr}}

\newcommand{\cen}{{\rm cen}}
\newcommand{\Pprob}{{\mathbb{P}}}

\newcommand{\dI}{{\mathbb{I}}}
\newcommand{\dR}{{\mathbb{R}}}

\newcommand{\rd}{{\rm d}}

\newcommand{\cov}{{\rm cov}}
\newcommand{\var}{{\rm var}}

\newcommand{\re}{{\rm e}}
\newcommand{\rA}{{\rm \scriptscriptstyle A}}
\newcommand{\rE}{{\rm \scriptscriptstyle E}}
\newcommand{\E}{{\mathbb{E}}}
\newcommand{\M}{{\mathbb{M}}}
\newcommand{\RR}{{\mathbb{R}}}
\newcommand{\NN}{{\mathbb{N}}}
\newcommand{\trans}[1]{{#1}^{\scriptscriptstyle{\rm T}}}
\newcommand{\inv}[1]{{#1}^{{-1}}}

\newcommand{\diag}{{\rm diag}}
\newcommand{\viz}{{{\it viz}.}}
\newcommand{\nit}{n_{\rm it}}
\newcommand{\GJ}{{\scriptscriptstyle\rm GJ}}
\newcommand{\GJOR}{{\scriptscriptstyle\rm GJOR}}
\newcommand{\idp}{p^{-1}\cdot\mathds{I}_v}
\newcommand{\barnit}{{\bar n}_{\rm it}}
\newcommand{\nohyphens}{\hyphenpenalty=10000\exhyphenpenalty=10000\relax}

\newcommand{\ie}{{\it i.e.\ }}
\newcommand{\eg}{{\it e.g.\ }}

\newcommand{\epsT}{{\epsilon}_{\scriptscriptstyle \rm T}}
\newcommand{\Mo}[1]{{{\rm Mo}(#1)}}
\newcommand{\hatMo}[1]{{{\rm \widehat{Mo}}(#1)}}
\newcommand{\kcut}{{k_{\rm th}}}
\newtheorem{prop}{Proposition}[section]
\newtheorem{corol}{Corollary}[section]
\newtheorem{remark}{Remark}[section]
\newtheorem{theo}{Theorem}[section]

\renewcommand{\qedsymbol}{$\boxdot$\ }

\def\lambdabar{\protect\@lambdabar}
\def\@lambdabar{%
\relax
\bgroup
\def\@tempa{\hbox{\raise.73\ht0
\hbox to0pt{\kern.25\wd0\vrule width.5\wd0
height.1pt depth.1pt\hss}\box0}}%
\mathchoice{\setbox0\hbox{$\displaystyle\lambda$}\@tempa}%
{\setbox0\hbox{$\textstyle\lambda$}\@tempa}%
{\setbox0\hbox{$\scriptstyle\lambda$}\@tempa}%
{\setbox0\hbox{$\scriptscriptstyle\lambda$}\@tempa}%
\egroup
}

\DeclareMathAlphabet\mathbfcal{OMS}{cmsy}{b}{n}
\DeclareMathAlphabet{\mathcalligra}{T1}{calligra}{m}{n}

\begin{document}
\maketitle

\begin{abstract}
In this paper we relate the matrix ${\frak S}_{\cal B}$ of the second moments of a 
spherically truncated normal multivariate to its full covariance
matrix $\Sigma$ and present an algorithm to invert the relation and
reconstruct $\Sigma$ from $\fS_\cB$. While the eigenvectors of $\Sigma$
are left invariant by the truncation, its eigenvalues are non--uniformly
damped. We show that the eigenvalues of $\Sigma$ can be reconstructed from their
truncated counterparts via a fixed point iteration, whose convergence we prove
analytically. The procedure requires the computation of
multidimensional Gaussian integrals over a Euclidean ball, for which we extend a
numerical technique, originally proposed by Ruben in 1962, based on a 
series expansion in chi--square distributions. In order to study the
feasibility of our approach, we examine the convergence rate of some
iterative schemes on suitably chosen ensembles of Wishart matrices. 
We finally discuss the practical difficulties arising in sample space
and outline a regularization of the problem based on perturbation theory. 
\end{abstract}

\section{Introduction}

It is more than forty years since Tallis \cite{tallis} worked out the moment--generating function of a
normal multivariate $X\equiv\{X_k\}_{k=1}^{v}\,\sim\,\cN_v(0,\Sigma)$, subject to
the conditional event
\begin{equation}
X\in\cE_v(\rho;\Sigma)\,,\qquad\qquad \cE_v(\rho;\Sigma) \equiv \,\{x\in\RR^v:\ \trans{x}\Sigma^{-1}x \le \rho\}\,.
\label{eq:ellipcond}
\end{equation}
The perfect match between the symmetries of the ellipsoid ${\cE}_v(\rho;\Sigma)$
and those of $\cN_v(0,\Sigma)$ allows for an exact analytic
result, from which the complete set of multivariate truncated moments can be
extracted upon differentiation. Consider for instance the matrix
$\fS_\cE(\rho;\Sigma)$ of the second truncated moments, expressing the
covariances among the components of $X$ within $\cE_v(\rho;\Sigma)$. From
Tallis' paper it turns out that 
\begin{equation}
\fS_\cE(\rho;\Sigma) = c_{\scriptscriptstyle\rm T}(\rho)\, \Sigma\,,\qquad
c_{\scriptscriptstyle\rm T}(\rho) \equiv\, \frac{F_{v+2}(\rho)}{F_v(\rho)}\,,
\label{eq:tallisres}
\end{equation}
with $F_v$ denoting the cumulative distribution function of a
$\chi^2$--variable with $v$ degrees of freedom. Inverting
eq.~(\ref{eq:tallisres}) -- so as to express $\Sigma$ as a function of
$\fS_\cE$ --  is trivial, since $c_{\scriptscriptstyle\rm T}(\rho)$ is a scalar damping
factor independent of $\Sigma$. In this paper, we shall refer to such inverse
relation as the {\it reconstruction} of $\Sigma$ from $\fS_\cE$. Unfortunately, life is not always so easy. In
general, the effects produced on the expectation of functions
of $X$ by cutting off the probability density outside
a generic domain $\cD\subset\RR^v$ can be hardly calculated in closed form,
especially if the boundary of $\cD$ is shaped in a way that breaks the
ellipsoidal symmetry of $\cN_v(0,\Sigma)$. Thus, for instance, unlike
eq.~(\ref{eq:tallisres}) the matrix of the second truncated moments 
is expected in general to display a non--linear/non--trivial dependence upon $\Sigma$.     

In the present paper, we consider the case where $\cD$ is a $v$--dimensional
Euclidean ball with center in the origin and square radius
$\rho$. Specifically, we discuss the reconstruction of $\Sigma$
from the matrix $\fS_\cB$ of the spherically truncated second moments. To this
aim, we need to mimic Tallis' calculation, with eq.~(\ref{eq:ellipcond}) 
replaced by the conditional event  
\begin{equation}
X\in\cB_v(\rho)\,,\qquad\qquad \cB_v(\rho) \equiv\, \{x\in\RR^v:\ \trans{x}x\le\rho\}\,.
\label{eq:sphercond}
\end{equation}
This is precisely an example of the situation described in the previous
paragraph: although $\cB_v(\rho)$ has a higher degree of symmetry than
$\cE_v(\rho;\Sigma)$, still there is no possibility of obtaining a closed--form
relation between $\Sigma$ and $\fS_\cB$, since $\cB_v(\rho)$ breaks the ellipsoidal
symmetry of $\cN_v(0,\Sigma)$: technically speaking, in this 
case we cannot perform any change of variable under the
defining integral of the moment--generating function, which may reduce the
dimensionality of the problem, as in Tallis' paper. 

In spite of that, the residual symmetries characterizing the truncated distribution help
simplify the  problem in the following respects: {\it i}) the reflection
invariance of the  whole set--up still yields $\E[X_k\,|\,X\in\cB_v(\rho)] = 0$
$\forall\,k$,  and  {\it ii}) the rotational invariance of $\cB_v(\rho)$
preserves the  possibility of defining the principal components of the
distribution just like in the unconstrained case.  
In particular, the latter property means that $\fS_\cB$ and $\Sigma$
share the same orthonormal eigenvectors. 
In fact, the reconstruction of $\Sigma$ from $\fS_\cB$ amounts to solving a
system  of non--linear  integral equations, having the eigenvalues $\lambda \equiv
\{\lambda_k\}_{k=1}^v$ of $\Sigma$ as unknown variables and the eigenvalues
$\mu\equiv\{\mu_k\}_{k=1}^v$ of $\fS_\cB$ as input parameters. In a lack
of analytic techniques to evaluate exactly the integrals involved, we 
resort to a numerical algorithm, of which we investigate feasibility,
performance and optimization.  

The paper is organized as follows. In sect.~2, we describe a couple
of examples illustrating the occurrence of spherical
truncations in practical situations. In sect.~3, we show that the
aforementioned integral equations have the 
analytic structure of a fixed point vector equation, that is to say
$\lambda = T(\lambda)$. This suggests to achieve the reconstruction of
$\lambda$ numerically via suitably chosen iterative schemes. In
sect.~4, we prove the convergence of the simplest of them by inductive
arguments, the validity of which relies upon the  monotonicity
properties of ratios of Gaussian integrals over $\cB_v(\rho)$. In
sect.~5, we review some numerical techniques for the computation of
Gaussian integrals over $\cB_v(\rho)$ with controlled systematic
error. These are based on and extend a classic work by
Ruben~\cite{ruben4} on the distribution of quadratic forms of
normal variables. For the sake of readability,  we defer proofs of
statements made in this sect. to Appendix A. In sect.~6, we report
on our numerical experiences: since the simplest iterative scheme,
namely the Gauss--Jacobi iteration, is too slow for practical
purposes, we investigate the performance of its improved version
based on over--relaxation; as expected, we find that the latter
has a higher convergence rate, yet it still slows down polynomially in
$1/\rho$ as $\rho\to 0$ and exponentially in $v$ as
$v\to\infty$; in order to reduce the slowing down, we propose an
acceleration technique, which boosts the higher components of the
eigenvalue spectrum. A series of Monte Carlo simulations enables
us to quantify the speedup. In sect.~7 we discuss the problems
arising when $\mu$ is affected by statistical uncertainty and propose a 
regularization technique based on perturbation theory. To conclude, 
we summarize our findings in sect.~8. 

\section{Motivating examples}

Spherical truncations of multinormal distributions may characterize
different kinds of experimental devices and may occur in 
various problems of statistical and convex analysis. In this
section, we discuss two motivating examples. 

\subsection{A two--dimensional {\it gedanken} experiment in Classical Particle Physics}

\begin{figure}
\centering
\includegraphics[width=0.6\textwidth]{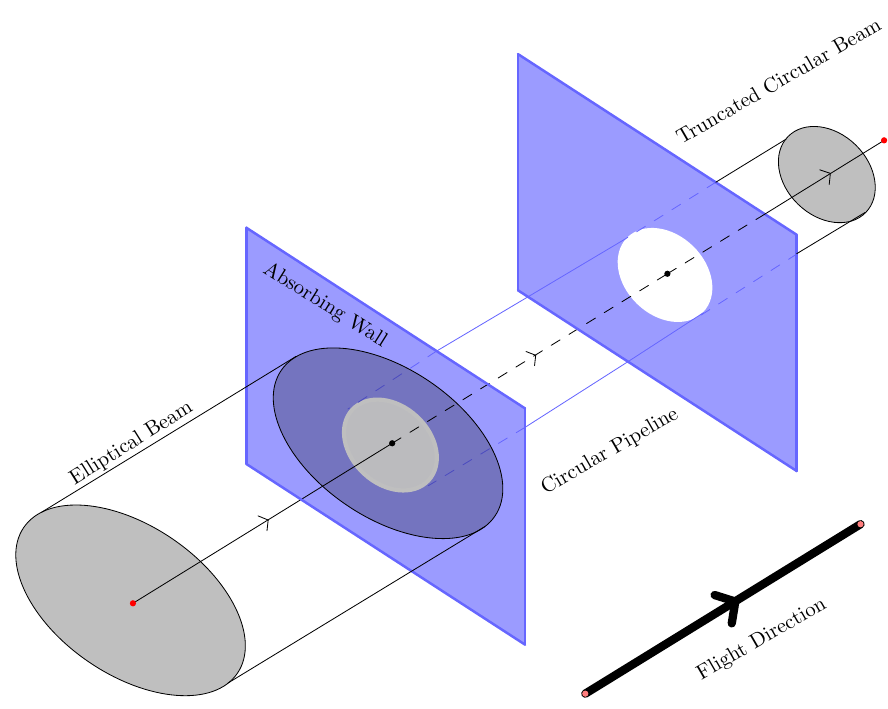}
\caption{\small A classical particle beam with elliptical
  transversal profile is cut off upon entering a circular coaxial pipeline.\label{fig:accel}}
\end{figure}

Consider the following ideal situation. An accelerator
physicist prepares an elliptical beam of classical particles with
Gaussian transversal profile. The experimenter knows {\it a priori} the
spatial distribution of the beam, \ie the covariance matrix
$\Sigma$ of the two--dimensional coordinates
of the particles on a plane orthogonal to their flight direction.
We can assume with no loss of generality that the transversal coordinate
system has origin at the maximum of the beam intensity and axes
along the principal components of the beam, thus it holds 
$\Sigma = \diag(\lambda_1,\lambda_2)$. The beam travels
straightforward until it enters a linear coaxial pipeline with
circular profile, schematically depicted in Fig.~\ref{fig:accel},
where the beam is longitudinally accelerated. While the outer part
of the beam is stopped by an absorbing wall, the inner part
propagates within the pipeline. At the end of the beam flight the
physicist wants to know if the transversal distribution of the
particles is changed, due to unknown disturbance factors arisen
within the pipeline. Accordingly, he measures again the spatial
covariance matrix of the beam. Unfortunately, the absorbing wall
has cut off the Gaussian tail, thus damping the covariance matrix
and making it not anymore directly comparable to the original
one. To perform such a comparison in the general case
$\lambda_1\ne \lambda_2$, the covariance matrix of the
truncated circular beam has to go through the reconstruction
procedure described in next sections.  

\subsection{A multivariate example: connections to Compositional Data Analysis}

Compositional Data Analysis (CoDA) has been the subject of a number of
papers, pioneered by J. Aitchison~\cite{aitchison} over the past
forty years. As a methodology of statistical investigation,
it finds application in all cases where the
main object of interest is a multivariate with
strictly positive continuous components to be
regarded as portions of a total amount $\kappa$ (the normalization
$\kappa=1$ is conventionally adopted in the mathematical
literature). In other words, compositional variates belong to the
$\kappa$--simplex 
\begin{equation}
\cS_v = \left\{z\in\dR^v_+: \quad |z|_1=\kappa\right\}\,,\qquad v\ge2\,,
\end{equation}
with $|z|_1 = \sum_{k=1}^v z_k$ the taxi--cab norm of
$z$, while compositions with different norms can be
always projected onto $\cS_v$ via the closure operator
$\cC\cdot x \equiv \{\kappa x_1/|x|_1,\ldots,\kappa
x_v/|x|_1\}$. There are countless types of compositional data, whose
analysis raises problems of interest for Statistics \cite{PawloTwo}, \eg geochemical data, balance sheet data, election
data, {\it etc}.

In order to measure distances on $\cS_v$, Aitchison introduced a
positive symmetric function \linebreak $d_\rA:\cS_v\times\cS_v\to\dR_+$, explicitly defined by
\begin{equation}
d_\rA(x,y) =
\sqrt{\frac{1}{2v}\sum_{i,k=1}^v\left[\log\left(\frac{x_i}{x_k}\right)
  - \log\left(\frac{y_i}{y_k}\right)\right]^2}\,.
\end{equation}
The Aitchison distance is a key tool in CoDA. It is scale
invariant in both its first and second argument, \ie it is left invariant by
redefinitions  $z\to\{\alpha z_1,\ldots,\alpha
z_v\}$ with $\alpha\in\dR_+$.  Accordingly, its support can be
extended to $\dR_+^v\times\dR_+^v$ by imposing 
\begin{equation}
d_\rA(x,y) \equiv d_\rA(\cC\cdot x,\cC\cdot y)\,,\qquad  x,y\in\dR_+^v\,.
\end{equation}
It was proved in \cite{isometric} that $d_\rA$ is
a norm--induced metric on $\cS_v$, provided the latter is given
an appropriate normed vector space structure.  Owing to the compositional constraint
$|\cdot |_1=\kappa$, it holds $\text{dim}\,\cS_v=v-1$. Accordingly, the
description of $\cS_v$ in terms of $v$ components is
redundant: an essential representation requires 
compositions to be properly mapped onto
$(v-1)$--tuples. Among various possibilities, the Isometric 
Logratio Transform (ILR) introduced in \cite{isometric}, is the 
only known map of this kind leaving $d_\rA$ invariant. More
precisely, the ILR fulfills
\begin{align}
d_\rA(x,y)  = d_\rE\left(\ilr(x),\ilr(y)\right)\,,\qquad d_\rE(u,v) \equiv \sqrt{\sum_{k=1}^{v-1}(u_k-v_k)^2}\,.
\label{eq:isometry}
\end{align}

It is known from \cite{AitShen} that if
$X\sim \log \cN_v(\mu,\Sigma)$ is a log--normal $v$--variate, then
$\cC\cdot X\sim L_v(\mu',\Sigma')$ is a logistic--normal
$v$--variate\footnote{the reader should notice that in \cite{AitShen} the simplex is
defined by $\cS_{v} = \{z\in\dR_+^v : |z|_1<1\}$, thus
  property 2.2 of \cite{AitShen} is here reformulated so as to take into
account such difference.}, with a known relation between $(\mu,\Sigma)$ and
$(\mu',\Sigma')$. Analogously, it is not difficult to show that
$\ilr\left(\cC\cdot X\right)\sim\cN_{v-1}(\mu'',\Sigma'')$ is a normal
$(v-1)$--variate, with $(\mu'',\Sigma'')$  related to
$(\mu',\Sigma')$ via the change of basis matrices derived
in~\cite{isometric}. Just to sum up, it holds
\begin{equation}
X\sim \log \cN_v(\mu,\Sigma)\quad\Rightarrow\quad \cC\cdot X\sim
L_v(\mu',\Sigma') \quad\Leftrightarrow\quad \ilr\left(\cC\cdot
  X\right)\sim\cN_{v-1}(\mu'',\Sigma'')\,. 
\label{eq:refimpl}
\end{equation}
Now, suppose that ${\it i})$ $X$ fulfills eq.~(\ref{eq:refimpl})
and has a natural interpretation as a composition, ${\it ii})$ a
representative set $\cD_X$ of observations of $X$ is given and ${\it iii})$
we wish to select from $\cD_X$ those units which are
compositionally closer to the centre of the distribution, according to the
Aitchison distance. To see that the problem is well posed,  we
first turn $\cD_X$ into a set $\cD_{\cC\cdot X} \equiv \{y: y=
\cC\cdot  x\ \text{and}\ x\in\cD_X\}$ of compositional
observations of $Y=\cC\cdot X$. Then, we consider the special
point $\cen[Y] = \cC\cdot\exp\{\E[\ln Y]\}$, representing the
centre of the distribution of $Y$ in a compositional sense:
$\cen[Y]$ minimizes the expression $\E[d^2_\rA(Y,\cen[Y])]$ over
$\cS_v$ and fulfills $\cen[Y] =
\ilr^{-1}(\E[\ilr(Y)])$, see ref.~\cite{PawloOne}. By virtue of eq.~(\ref{eq:refimpl}) this
entails $\ilr(\cen[Y]) = \E[\ilr(Y)] = \mu''$. In order to select
the observations which are closer to $\cen[Y]$, we set a threshold
$\delta>0$ and consider only those elements $y\in\cD_{\cC\cdot X}$
fulfilling $d_\rA^2(y,\cen[y])<\delta$, with $\cen[y]$ being a
sample estimator of $\cen[Y]$ on $\cD_{\cC\cdot X}$. Such
selection rule operates a well--defined truncation of the
distribution of $Y$.
Moreover, in view of eqs.~(\ref{eq:isometry}) and (\ref{eq:refimpl}), we have
\begin{equation}
\Pprob\left[d^2_\rA(Y,\cen[Y])<\delta\ |\ Y\sim
  L_v(\mu',\Sigma')\right] =
\Pprob\left[d^2_\rE(Z,\mu'')<\delta\ |\ Z\sim \cN_{v-1}(\mu'',\Sigma'')\right]\,,
\end{equation}
with $Z=\ilr(\cC\cdot X)$. As a consequence, we see that a
compositional selection rule based on the Aitchison distance and
eq.~(\ref{eq:refimpl}) is equivalent to a spherical truncation of a multinormal
distribution. Obviously, once $Z$ has been
spherically truncated, the covariance matrix of the remaining data
is damped, thus an estimate of the full covariance matrix requires
the reconstruction procedure described in next sections. 

\subsection{General covariance reconstruction problem}

The examples discussed in the previous subsections are special cases of a more general inverse problem, 
namely the reconstruction of the covariance matrix $\Sigma$ of a normal multivariate 
$X$ on the basis of the covariance matrix ${\frak S}_\cD$ of $X$ truncated to some (convex) region $\cD$. 
This is the simplest yet non--trivial inverse problem, which can be naturally associated to the normal 
distribution. The case $\cD=\cB_v(\rho)$ corresponds to a set--up where theoretical and practical 
aspects of the problem can be investigated with relatively modest mathematical effort. It is certainly
a well--defined framework where to study regularization techniques for non--linear inverse problems in 
Statistics, for which there is still much room for interesting work \cite{englreg,cavalier}. 

\section{Definitions and set--up}

Let $X\in\RR^v$ be a random vector with jointly normal distribution
$\cN_v(0,\Sigma)$ in $v\ge 1$ dimensions. The probability that $X$ falls
within $\cB_v(\rho)$ is measured by the Gaussian integral 
\begin{equation}
\alpha(\rho;\Sigma)\, \equiv\, \Pprob\left[X\in \cB_v(\rho)\right]\, =
 \, \frac{1}{(2\pi)^{v/2}|\Sigma|^{1/2}}\int_{\cB_v(\rho)}\rd^vx\ \exp\left\{-\frac{1}{2} \trans{x}\inv{\Sigma}x\right\}\,.
\label{eq:alpha1}
\end{equation}
Since $\Sigma$ is symmetric positive definite, it has orthonormal
eigenvectors $\Sigma v^{(i)} = \lambda_i v^{(i)}$. Let us denote by
$R\equiv\{v^{(j)}_i\}_{i,j=1}^{v}$ the orthogonal matrix having
these vectors as columns and by $\Lambda \equiv \diag(\lambda)=\trans{R}\Sigma
R$ the diagonal counterpart of $\Sigma$. From the invariance of $\cB_v(\rho)$
under rotations, it follows that $\alpha$ depends upon $\Sigma$ just by way of
$\lambda$. Accordingly, we rename the Gaussian probability content of
$\cB_v(\rho)$ as 
\begin{equation}
\alpha(\rho;\lambda) \equiv \int_{\cB_v(\rho)}\rd^v x \,\prod_{m=1}^v
\delta(x_m,\lambda_m)\,,\qquad\qquad \delta(y,\eta)\, =\, \frac{1}{\sqrt{2\pi\eta}}\exp\left\{-\frac{y^2}{2\eta}\right\}\,.
\label{eq:alpha2}
\end{equation}
Note that eq.~(\ref{eq:alpha2}) is not sufficient to fully characterize the random
vector $X$ under the spherical constraint, for which we need to calculate the
distribution law $\Pprob[X\in A|X\in\cB_v(\rho)]$ as a function of
$A\subset\RR^v$. Alternatively, we can describe $X$ in terms of the complete
set of its truncated moments 
\begin{equation}
m_{k_1\dots k_v}(\rho;\Sigma)\, \equiv\, \E[X_1^{k_1}\dots
X_v^{k_v}\,|\,X\in\cB_v(\rho)]\,,\qquad \{k_i\}_{i=1,\ldots,v}\in \NN^v\,,
\end{equation}
As usual, these can be all obtained from the moment--generating function
\begin{equation}
\alpha\, m(t) =\, \frac{1}{(2\pi)^{v/2}|\Sigma|^{1/2}}\int_{\cB_v(\rho)}\rd^v x\,
 \exp\left\{\trans{t}x-\frac{1}{2} \trans{x}\inv{\Sigma}x\right\}\,,\qquad t\in \RR^v\,,
\end{equation}
by differentiating the latter an arbitrary number of times with respect to the components of $t$, \viz
\begin{equation}
m_{k_1\dots k_v}(\rho;\Sigma) =\, \frac{\partial^{\,k_1+\ldots+k_v}\, m(t)}{(\partial t_1)^{k_1}\dots (\partial t_v)^{k_v}}\biggr|_{t=0}\,.
\end{equation}
It will be observed that $m(t)$ is in general not invariant under
rotations of $t$. Therefore, unlike $\alpha$, the moments $m_{k_1\ldots k_v}$ depend
effectively on both $\lambda$ and $R$. For instance, for the matrix of the second moments
  $\fS_\cB \equiv \{\partial^2 m/\partial t_i\partial t_j|_{t=0}\}_{i,j=1}^v$ such dependence amounts to
\begin{equation}
\alpha\, ({\fS_\cB})_{ij} \,=\, \sum_{k,\ell = 1}^v R_{ki}R_{\ell j} \int_{\cB_v(\rho)} \rd^v x\ x_kx_\ell\,\prod_{m=1}^v \delta(x_m,\lambda_m)\,.
\label{eq:truncsecondmom}
\end{equation}
By parity, the only non--vanishing terms in the above sum are those
with $k=\ell$. Hence, it follows that $\Sigma$ and $\fS_\cB$ share $R$ as a
common diagonalizing matrix. In other words, if $M\equiv\diag(\mu)$ is the diagonal
matrix of the eigenvalues of $\fS_\cB$, then $M =
\trans{R}\fS_\cB R$. Moreover, $\mu_k$ is related to $\lambda_k$ by 
\begin{equation}
\mu_k = \lambda_k\frac{\alpha_k}{\alpha}\,,\qquad\qquad \alpha_k(\rho;\lambda)
\equiv \int_{\cB_v(\rho)}\rd^v x\
\frac{x_k^2}{\lambda_k}\,\prod_{m=1}^v\delta(x_m,\lambda_m)\,,\qquad k = 1,\ldots,v\,.
\label{eq:truncspectrum}
\end{equation}
The ratios $\alpha_k/\alpha$ are naturally interpreted as adimensional correction
factors to the eigenvalues of $\Sigma$, so they play the same r\^ole as 
$c_{\scriptscriptstyle \rm T}(\rho)$ in eq.~(\ref{eq:tallisres}). However,
$\alpha_k/\alpha$ depends explicitly on the subscript~$k$, thus each eigenvalue
is damped differently from the others as a consequence of the condition $X\in\cB_v(\rho)$.   
\vskip 0.4cm
\begin{remark}
In practical terms, eqs.~(\ref{eq:truncsecondmom})--(\ref{eq:truncspectrum})
tell us that estimating the sample covariance matrix of $X\sim\cN_v(0,\Sigma)$ from a spherically
truncated population $\{x^{(m)}\}_{m=1}^M$, made of $M$
independent units, via the classical estimator $Q_{ij} = (M-1)^{-1}\sum_{m=1}^M(x_i^{(m)}-{\tilde
  x}^{(m)}_i)(x_j^{(m)}-{\tilde x}^{(m)}_j)$, being $\tilde x = M^{-1}\sum_{m=1}^M
x^{(m)}$ the sample mean, yields a damped result. Nonetheless, the
damping affects only the eigenvalues of the estimator, whereas its eigenvectors are
left invariant.  
\end{remark}

\subsection{Montonicity properties of ratios of Gaussian integrals}

Eqs.~(\ref{eq:alpha2}) and (\ref{eq:truncspectrum}) suggest to introduce a
general notation for the Gaussian integrals over $\cB_v(\rho)$, under the
assumption $\Sigma = \Lambda$. So, we define
\begin{equation}
\alpha_{k\ell m\dots}(\rho;\lambda) \equiv \int_{\cB_v(\rho)}\rd^v x\ \frac{x_k^2}{\lambda_k}\,\frac{x_\ell^2}{\lambda_\ell}\,\frac{x_m^2}{\lambda_m}\dots\,\prod_{n=1}^v\delta(x_n,\lambda_n)\,,
\end{equation}
with each subscript $q$ on the {\it l.h.s.} addressing an additional factor of $x_q^2/\lambda_q$
under the integral sign on the {\it r.h.s.} (no subscripts means
$\alpha$). Several analytic properties of such integrals are
discussed in ref.~\cite{palovar}. Here, we lay emphasis on some issues concerning
the monotonicity trends of the ratios $\alpha_k/\alpha$. Specifically, 
\begin{prop}[monotonicities]
\label{prop:mon}
Let $\lambda_{(k)}\equiv\{\lambda_i\}_{i=1,\ldots,v}^{i\ne k}$ denote the set of the
full eigenvalues without $\lambda_k$. The ratios $\alpha_k/\alpha$ fulfill the
following properties:  
\begin{itemize}
\item[$p_1${\rm )}]{\ $\lambda_k\dfrac{\alpha_k}{\alpha}(\rho;\lambda)$ is a
    monotonic increasing function of $\lambda_k$ at fixed $\rho$ and $\lambda_{(k)}$;}
\item[$p_2${\rm )}]{\ $\dfrac{\alpha_k}{\alpha}(\rho;\lambda)$ is a monotonic
    decreasing function of $\lambda_k$ at fixed $\rho$ and $\lambda_{(k)}$;}
\item[$p_3${\rm )}]{\ $\dfrac{\alpha_k}{\alpha}(\rho;\lambda)$ is a monotonic
    decreasing function of $\lambda_i$ $(i\ne k)$ at fixed $\rho$ and $\lambda_{(i)}$,}
\end{itemize}
where an innocuous abuse of notation has been made on writing
$\frac{\alpha_k}{\alpha}(\rho;\lambda)$ in place of $\alpha_k(\rho;\lambda)/\alpha(\rho;\lambda)$.
\end{prop}
\begin{proof}
Let the symbol $\partial_k \equiv \partial/\partial\lambda_k$ denote a
derivative with respect to $\lambda_k$. In order to prove property $p_1$),
we apply the chain rule of differentiation to $\lambda_k\alpha_k/\alpha$ and then
we pass $\partial_k$ under the integral
sign in $\partial_k\alpha$ and $\partial_k\alpha_k$. In this way, we obtain
\begin{align}
\partial_k\left(\lambda_k\frac{\alpha_k}{\alpha}\right) & = \
\frac{1}{2}\left(\frac{\alpha_{kk}}{\alpha}-\frac{\alpha_k^2}{\alpha^2}\right)
= \frac{1}{2\lambda_k^2}\left\{\E[X_k^4\,|\,X\in \cB_v(\rho)] - 
  \E[X_k^2\,|\,X\in \cB_v(\rho)]^2\right\} \nonumber\\[2.0ex] 
& = \
\frac{1}{2\lambda_k^2}\var\left(X_k^2\,|\,X\in\cB_v(\rho)\right) \ge 0\,.
\label{eq:proofp1}
\end{align}
Moreover, since the truncated marginal density of $X_k^2$ is positive within a set of
non--zero measure in $\RR$, the monotonic trend of $\lambda_k\alpha_k/\alpha$
in $\lambda_k$ is strict. \qedsymbol Properties $p_2$)
and $p_3$) are less trivial than $p_1$). Indeed, the same reasoning as above now
yields on the one hand 
\begin{align}
\lambda_k\partial_k\left(\frac{\alpha_k}{\alpha}\right)
& = \partial_k\left(\lambda_k\frac{\alpha_k}{\alpha}\right) -
\frac{\alpha_k}{\alpha} \nonumber\\[2.0ex]
& = \frac{1}{2\lambda_k^2}\left\{\var\left(X_k^2\,|\,X\in\cB_v(\rho)\right)
- 2\lambda_k\E[X_k^2\,|\,X\in\cB_v(\rho)]\right\}\le 0\,,
\label{eq:varineq}
\end{align}
and on the other
\begin{equation}
\hskip 0.1cm\lambda_i\partial_i\left(\frac{\alpha_k}{\alpha}\right) =
\frac{1}{2}\left(\frac{\alpha_{ik}}{\alpha}-
  \frac{\alpha_i\alpha_k}{\alpha^2}\right) =
\frac{1}{2\lambda_i\lambda_k}\cov\left(X_i^2,X_k^2\,|\,X\in\cB_v(\rho)\right)\le 0\,\qquad
(i\ne k)\,.
\label{eq:covineq}
\end{equation}
Though not {\it a priori} evident, the {\it r.h.s.} of both
eqs.~(\ref{eq:varineq}) and (\ref{eq:covineq})
is negative (and vanishes in the limit $\rho\to\infty$). The
inequalities $\var(X_k^2)\le 2\lambda_k\E[X_k^2]$ within Euclidean
balls have been first discussed in \cite{palovar}, while the
inequalities $\cov(X_j^2,X_k^2)$ within generalized Orlicz balls
have been discussed in refs.~\cite{ball:2003,orlicz1} for the case where the probability 
distribution of $X$ is flat instead of being normal. More
recently, a complete proof of both inequalities has been given in
\cite{mukerjee}. Despite the technical difficulties in proving
them, their meaning should be intuitively clear. The 
variance inequality quantifies the squeezing affecting $X_k^2$ as
a consequence of the truncation (in the unconstrained case it
would be $\var(X_k^2)=2\lambda_k^2$). The covariance inequality
follows from the opposition arising among  the square components
in proximity of the boundary of $\cB_v(\rho)$. Indeed, if
$X_j^2\nearrow\rho$,  then $X^2_k\searrow 0$ $\forall\, k\ne j$ in
order for $X$ to stay within $\cB_v(\rho)$. 
\end{proof}

\subsection{Definition domain of the reconstruction problem}

\noindent A consequence of Proposition~\ref{prop:mon} is represented by the following
\begin{corol}
Given $v$, $\rho$ and $\lambda$, $\mu_k$ is bounded by
\begin{equation}
\frac{\rho}{r\left(v,\frac{\rho}{2\lambda_k}\right)}\le\,\mu_k\,\le \frac{\rho}{3}\,,\qquad\qquad
r(v,z) \equiv (2v+1)\dfrac{M\left(v,v+{1}/{2},z\right)}{M\left(v,v+{3}/{2},z\right)}\,,
\label{eq:mubounds1}
\end{equation}
with $M$ denoting the Kummer function, {\it viz.}
\begin{equation}
M(a,b,z) = \sum_{n=0}^\infty \frac{1}{n!}\frac{(a)_n}{(b)_n}\,z^n\,,\qquad
(x)_n \equiv \frac{\Gamma(x+n)}{\Gamma(x)}\,.
\end{equation}
\end{corol}
\begin{proof}
The upper bound of eq.~(\ref{eq:mubounds1}) corresponds to the value of
$\mu_k$ in the $v$--tuple limit $\lambda_k\to\infty$, $\lambda_{(k)}\to
\{0\,\ldots,0\}$. This is indeed the maximum possible value allowed for
$\mu_k$ according to properties $p_1$) and $p_3$) of
Proposition~\ref{prop:mon}. In order to perform this limit, we  observe that 
\begin{equation}
\lim_{\eta\to 0^+}\delta(y,\eta) = \delta(y)\,,
\end{equation}
with the $\delta$ symbol on the {\it r.h.s.} representing the Dirac delta function
(the reader who is not familiar with the theory of distributions may refer for
instance to ref.~\cite{joshi} for an introduction).  Accordingly,
\begin{align}
\lim_{\lambda_k\to\infty}\, \lim_{\lambda_{(k)}\to\{0,\ldots,0\}}\,\mu_k \ =\ 
{\displaystyle{\int_{-\sqrt{\rho}}^{+\sqrt{\rho}}\rd
  x_k\ x_k^2}}\biggl/{\displaystyle{\int_{-\sqrt{\rho}}^{+\sqrt{\rho}}\rd x_k}} = \frac{\rho}{3}\,.
\end{align}
The lower bound corresponds instead to the value taken by $\mu_k$ as
$\lambda_{(k)}\to \{\infty\,\ldots,\infty\}$ and $\lambda_k$ is kept fixed. In this limit, all the Gaussian
factors in the probability density function except the $k^{\rm th}$ one
flatten to one. Hence, 
\begin{align}
& \lim_{\lambda_{(k)}\to \{\infty\,\ldots,\infty\}}\mu_k \ =
 \lim_{\lambda_{(k)}\to
  \{\infty\,\ldots,\infty\}}\frac{\displaystyle{\int_{-\sqrt{\rho}}^{+\sqrt{\rho}}\rd
    x_k\,x_k^2\,\delta(x_k,\lambda_k)\
    \alpha^{(v-1)}\left(\rho-x_k^2;\lambda_{(k)}\right)}}{\displaystyle{\int_{-\sqrt{\rho}}^{+\sqrt{\rho}}\rd
    x_k\ \delta(x_k,\lambda_k)\
    \alpha^{(v-1)}\left(\rho-x_k^2;\lambda_{(k)}\right) }} \nonumber
\end{align}
\begin{align}
& \hskip 1.0cm = \frac{\displaystyle{\int_{-\sqrt{\rho}}^{+\sqrt{\rho}}\rd
    x_k\, x_k^2\,\re^{-\frac{x_k^2}{2\lambda_k}}\,\left(\rho-x_k^2\right)^{v-1}}}{\displaystyle{\int_{-\sqrt{\rho}}^{+\sqrt{\rho}}\rd
    x_k\
    \re^{-\frac{x_k^2}{2\lambda_k}}\,\left(\rho-x_k^2\right)^{v-1}}}
= \rho \frac{\displaystyle{\int_{0}^{1}\rd
    x_k\,x_k^2\,\re^{-\frac{\rho}{2\lambda_k}x_k^2}\,\left(1-x_k^2\right)^{v-1}}}{\displaystyle{\int_{0}^{1}\rd
    x_k\
    \re^{-\frac{\rho}{2\lambda_k}x_k^2}\,\left(1-x_k^2\right)^{v-1}}}\,.
\end{align}
Numerator and denominator of the rightmost ratio are easily recognized to be
integral representations of Kummer functions (see {\it e.g.} ref.~\cite[ch.~13]{abramowitz}).
\end{proof}
\noindent The upper bound of eq.~(\ref{eq:mubounds1}) can be sharpened, as clarified by
the following
\begin{prop}[Bounds on the truncated moments] Let $v$, $\rho$ and $\lambda$ be
  given. If $\{i_1,\ldots,i_v\}$ is a permutation of $\{1,\ldots,v\}$ such that 
 $\mu_{i_1}\le\mu_{i_2}\le\ldots\le\mu_{i_v}$, then the following upper bounds hold:
\begin{equation}
i)\quad \sum_{k=1}^v \mu_k \ \le\ \rho\,;\qquad\qquad ii)\quad \mu_{i_k} \ \le\
\frac{\rho}{v-k+1}\,,\quad k=1,\ldots,v\,.
\label{eq:mubounds2}
\end{equation}
\end{prop}
\begin{proof}
The overall upper bound on the sum of truncated moments follows from
\begin{equation}
\sum_{k=1}^v \mu_k \ =\  \frac{1}{\alpha}\sum_{k=1}^v \lambda_k\alpha_k \ = \
\frac{1}{\alpha}\int_{\cB_v(\rho)}\rd^v x\ \left(\sum_{k=1}^v x_k^2\right)\,\prod_{m=1}^v
\delta(x_m,\lambda_m)\,\le \rho\,.
\label{eq:muboundsest1}
\end{equation}
At the same time, the sum can be split and bounded from below by
\begin{equation}
\sum_{k=1}^v \mu_{i_k}\  =\  \sum_{k=1}^n\mu_{i_k} +
\sum_{k=n+1}^v\mu_{i_k}\ \ge\ \sum_{k=1}^n\mu_{i_k}
+(v-n)\mu_{i_{n+1}}\,,\qquad n = 0,1,\ldots,v-1\,.
\label{eq:muboundsest2}
\end{equation}
The single upper bounds on the $\mu_k$'s are then obtained from eqs.~(\ref{eq:muboundsest1})%
--(\ref{eq:muboundsest2}). It will be noted that eq.~(\ref{eq:mubounds2}) {\it ii}) is sharper
than the upper bound of eq.~(\ref{eq:mubounds1}) only for $v>3$ and $k<v-2$.
\end{proof}

\begin{figure}[t!]
  \begin{minipage}[!t]{0.49\textwidth}
    \hskip 0.5cm
    \includegraphics[width=0.8\textwidth]{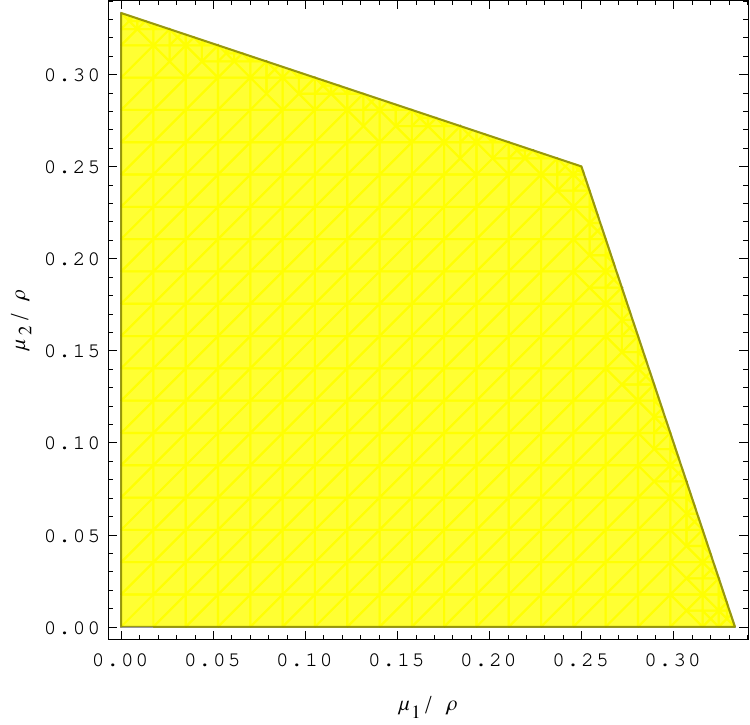}
    \hfill
  \end{minipage}
  \hskip 1.0cm\begin{minipage}[!t]{0.49\textwidth}
    \vskip -1.0cm
    \hfill
    \includegraphics[width=0.8\textwidth]{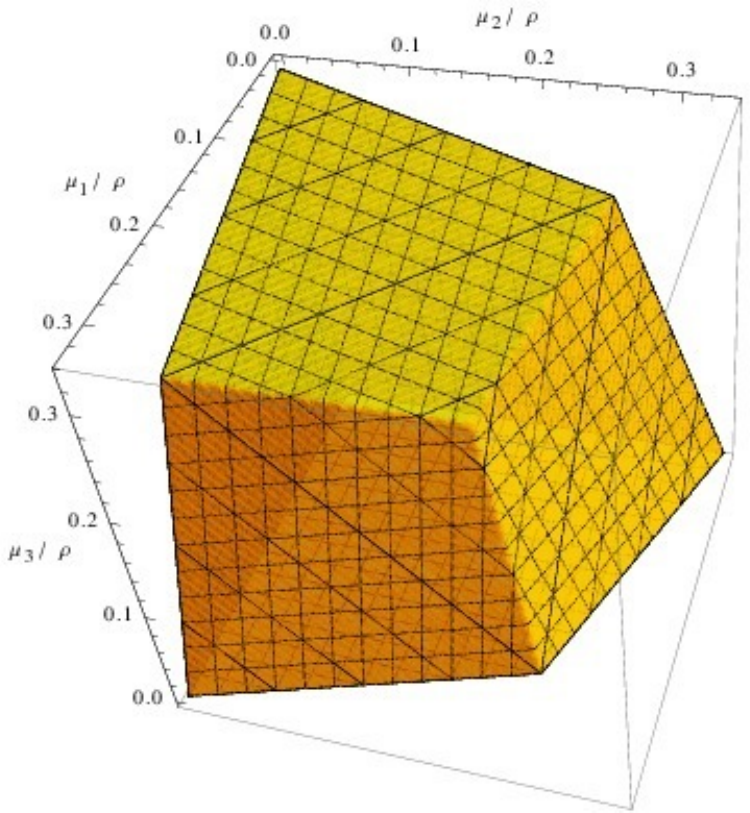}
    \hbox{\hskip 1.0cm \phantom{a}}
  \end{minipage}
  \vskip -0.2cm
  \caption{ \footnotesize {\it Left}: Numerical reconstruction of
    $\cD(\tau_\rho^{-1})$ in $v=2$ dimensions. {\it Right}: Numerical
    reconstruction of $\cD(\tau_\rho^{-1})$ in $v=3$ dimensions.}
  \label{fig:mubounds}
\end{figure}

From now on, we shall assume with no loss of generality, that the eigenvalues
of $\Sigma$ are increasingly  ordered, namely $0<\lambda_1\le\dots\le
\lambda_v$ (we can always permute the labels of the
coordinate axes, so as to let this be the case). An important aspect related
to the eigenvalue ordering is provided by the following 
\begin{prop}[Eigenvalue ordering]
Let $v$, $\rho$ and $\lambda$ be given. If $\lambda_{1}\le
\lambda_{2}\le\ldots\le\lambda_{v}$, then it holds $\mu_{1}\le
\mu_{2}\le\ldots\le\mu_{v}$ as well. 
\end{prop}
\begin{proof}
In order to show that the spherical truncation does not violate the eigenvalue ordering, we
make repeated use of the monotonicity properties of Proposition~3.1. Specifically, if
$i<j$ then 
\begin{align}
\mu_i & \ = \,
\lambda_i\frac{\alpha_i}{\alpha}(\rho;\{\lambda_1,\ldots,\lambda_i,\ldots,\lambda_j,\ldots,\lambda_v\})
\nonumber\\[2.0ex]
\phantom{\mu_1}
& \ \le\,
\lambda_j\frac{\alpha_i}{\alpha}(\rho;\{\lambda_1,\ldots,\lambda_j,\ldots,\lambda_j,\ldots,\lambda_v\})
\quad\quad\hskip0.2ex \color{Black}{\leftrightsquigarrow}\quad\text{increasing
  monotonicity of } \lambda_i\frac{\alpha_i}{\alpha} \nonumber\\[2.0ex]
\phantom{\mu_i}
& \ =\,
\lambda_j\frac{\alpha_j}{\alpha}(\rho;\{\lambda_1,\ldots,\lambda_j,\ldots,\lambda_j,\ldots,\lambda_v\})
\quad\quad \color{Black}{\leftrightsquigarrow}\quad\text{exchange
  symmetry } i\leftrightarrow j
\nonumber
\end{align}
\begin{align}
& \ \le\,
\lambda_j\frac{\alpha_j}{\alpha}(\rho;\{\lambda_1,\ldots,\lambda_i,\ldots,\lambda_j,\ldots,\lambda_v\})
\quad\quad\hskip0.11ex \color{Black}{\leftrightsquigarrow}\quad\text{decreasing
 monotonicity of } \frac{\alpha_j}{\alpha}
\nonumber\\[2.0ex]
& \ =\, \mu_j\,,
\end{align}
where the symbol ``$\leftrightsquigarrow$'' is used to explain where the
inequality sign preceding it comes from, and the ``exchange symmetry'' refers to the formal
property of the one--index Gaussian integrals over $\cB_v(\rho)$ to fulfill 
$\alpha_i(\rho;\{\lambda_1,\ldots,\lambda_i,\ldots,\lambda_j,\ldots,\lambda_v\}) =
\alpha_j(\rho;\{\lambda_1,\ldots,\lambda_j,\ldots,\lambda_i,\ldots,\lambda_v\})$.
\end{proof}
Let us now focus on eqs.~(\ref{eq:truncspectrum}). They have to be
solved in order to reconstruct $\lambda$ from $\mu$.  Formally, if we
introduce a family of truncation operators $\tau_\rho:\RR^v_+\to\RR^v_+$ (parametrically
depending on $\rho$), such that
\begin{align}
(\tau_\rho\cdot\lambda)_k \ \equiv \
  \lambda_k\frac{\alpha_k}{\alpha}(\rho;\lambda)\,,\qquad k=1,\dots,v\,,
\end{align}
then the reconstruction of $\lambda$ from $\mu$ amounts to calculating $\lambda =
\tau_\rho^{-1}\cdot\mu$. One should be aware that $\tau_\rho$ is not a
surjective operator in view of Corollary~3.1 and Proposition~3.2. Therefore,
$\tau_\rho^{-1}$ is only defined within a bounded domain
$\cD(\tau_\rho^{-1})$. If we define
\begin{align}
\cD_0 & = \left\{\mu\in\dR^v_+:\quad \mu_1\le\ldots\le\mu_v \ \text{
    and } \right. \nonumber\\[1.0ex]
& \hskip 2.55cm \left. \mu_k = \lambda_k \frac{\alpha_k}{\alpha}
\ \text{ for } \  k = 1,\ldots,v \ \text{ and for some } \lambda\in\dR^v_+\right\}\,,
\end{align}
then we have $\cD(\tau_\rho^{-1}) = \{\mu: \mu = \sigma\cdot \mu_0 \text{ for
some } \mu_0\in\cD_0 \text{ and } \sigma\in S_v \}$, where $S_v$
is the set of permutations of $v$ elements. From Proposition 3.2
we conclude $\cD_0\subseteq H_v(\rho)$, being
\begin{equation}
H_v(\rho) \ \equiv \ \left\{x\in \RR^v_+:\
  x_k\,\le\,\min\left\{\frac{\rho}{3},\frac{\rho}{v-k+1}\right\}\ \text{ and
  }\ \sum_{k=1}^vx_k\le\rho\,,\quad \forall\ k\right\}\,.
\end{equation}
In fact, there are vectors $\mu\in\RR^v_+$ fulfilling
$\mu\in H_v(\rho)$ and $\mu\notin \cD_0$, thus we conclude that
$\cD_0$ is a proper subset of $H_v(\rho)$. Numerical experiences based on the 
techniques discussed in the next sections show indeed that
\begin{equation}
\cD(\tau_\rho^{-1}) = \bigcap_{k=1}^v\biggl\{\,\mu\in\RR^v_+:\quad \sum_{j\ne
  k}^{1\ldots v}\mu_j + 3\mu_k \le \rho\,\biggr\}\,.
\label{eq:domain}
\end{equation}
A graphical representation of eq.~(\ref{eq:domain}) in $v=2$ and $v=3$ dimensions is depicted in 
Fig.~\ref{fig:mubounds}. The reader should note that until sect.~7 we shall always 
assume that $\mu$ comes from the application of $\tau_\rho$ to some $\lambda$, thus 
$\mu\in\cD(\tau_\rho^{-1})$ by construction. 

Now, we observe that eqs.~(\ref{eq:truncspectrum}) can be written
in the equivalent form
\begin{align}
& \lambda = T(\lambda;\mu;\rho)\,,\label{eq:fixedpoint0}\\[0.5ex]
& T:\ \RR^v_+\times\RR^v_+\times \RR_+\,\to\,\RR^v_+\,;\qquad  T_k(\lambda;\mu;\rho) = \mu_k\frac{\alpha}{\alpha_k}(\rho;\lambda)\,,\qquad k=1,\dots,v\,.
\end{align}
Since $\rho$ and $\mu$ are (non--independent) input parameters for the
covariance reconstruction problem (and in order to keep the notation light),
in the sequel we shall leave the dependence of $T$ upon $\rho$ and $\mu$ implicitly
understood, i.e., we shall write eq.~(\ref{eq:fixedpoint0}) as $\lambda =
T(\lambda)$. Hence, we see that the full eigenvalue spectrum $\lambda$ is a fixed
point of the operator $T$. This suggests to obtain it as the limit of a sequence  
\begin{align}
\label{eq:fixedpoint}
&\lambda^{(0)} =\, \mu\,,\qquad \lambda^{(n+1)}\, = \, T(\lambda^{(n)})\,,\qquad n = 0,1,\dots\,,\\[2.0ex]
& \lambda \,=\, \lim_{n\to\infty} \lambda^{(n)}\,,
\label{eq:iterlaw}
\end{align}
provided this can be shown to converge. Note that since $\alpha_k<\alpha$, it follows
that $T_k(\lambda^{(n)})>\mu_k$ $\forall n$, so the sequence is bounded from
below by $\mu$. In particular, this holds for $n=0$. Therefore, the sequence
moves to the right direction at least at the beginning. A formal proof of
convergence, based on the monotonicity properties stated by Proposition~3.1,
is given in next section.    

\section{Convergence of the fixed point equation}

We split our argument into three propositions, describing different
properties of the sequence $\lambda^{(n)}$. They assert respectively that {\it i})
the sequence is component--wise monotonic increasing; {\it ii}) the sequence
is component--wise bounded from above by any fixed point of $T$; {\it iii}) if
$T$ has a fixed point, this must be unique. Statements {\it i}) and {\it ii})
are sufficient to guarantee the convergence of the sequence to a finite
limit (the unconstrained spectrum is a fixed point of $T$).  In addition,
the limit is easily recognized to be a fixed point of $T$. Hence, statement
{\it iii}) guarantees that the sequence converges to the unconstrained
eigenvalue spectrum. We remark that all the monotonicities discussed in
Proposition 3.1 are strict, i.e. the ratios $\alpha_k/\alpha$ have no 
stationary points at finite $\rho$ and $\lambda$, which is crucial for the 
proof.  

\begin{prop}[Increasing monotonicity] Given $v$, $\rho$ and $\mu\in
  \cD(\tau_\rho^{-1})$, the sequence $\lambda^{(0)}=\mu$,
  $\lambda^{(n+1)} = T(\lambda^{(n)})$, $n=0,1,\ldots$ is monotonic increasing,
  viz. $\lambda_k^{(n+1)}>\lambda_k^{(n)}$ \ $\forall\, k=1,\ldots,v$. 
\end{prop}
\begin{proof}
The proof is by induction. We first notice that
\begin{equation}
\lambda_k^{(1)} = T_k(\lambda^{(0)}) = T_k(\mu) =
\mu_k\frac{\alpha}{\alpha_k}(\rho;\mu)> \mu_k = \lambda_k^{(0)}\,, \qquad k=1,\ldots,v\,,
\end{equation}
the inequality following from $\alpha_k(\rho;\mu)<\alpha(\rho;\mu)$. Suppose now that the
property of increasing monotonicity has been checked off up to the $n^{\rm th}$
element of the sequence. Then, 
\begin{equation}
\lambda_k^{(n+1)} = \mu_k\frac{\alpha}{\alpha_k}(\rho;\lambda^{(n)}) >
\mu_k\frac{\alpha}{\alpha_k}(\rho;\lambda^{(n-1)}) = \lambda_k^{(n)}\,,
\end{equation}
the inequality following this time from the inductive hypothesis and from
properties $p_2$) and $p_3$) of Proposition~3.1.
\end{proof}

\begin{prop}[Boundedness] Given $v$, $\rho$ and $\mu\in
  \cD(\tau_\rho^{-1})$, the sequence $\lambda^{(0)}=\mu$,
  $\lambda^{(n+1)} = T(\lambda^{(n)})$, $n=0,1,\ldots$ is bounded from above,
  {\it viz.} $\lambda_k^{(n)}< \lambda_k^{*}$ \ $\forall\, k=1,\ldots,v$, being
  $\lambda^*$ a fixed point of $T$.
\end{prop}
\begin{proof}
We proceed again by induction. We first notice that
\begin{equation}
\lambda_k^{(0)} = \mu_k <  \mu_k\frac{\alpha}{\alpha_k}(\rho;\lambda^*) = \lambda_k^{*}\,, \qquad k=1,\ldots,v\,,
\end{equation}
the inequality following as previously from $\alpha_k(\rho;\lambda^*)<\alpha(\rho;\lambda^*)$. Suppose now that the
property of boundedness has been checked off up to the $n^{\rm th}$
element of the sequence. Then, 
\begin{equation}
\lambda_k^{(n+1)} = \mu_k\frac{\alpha}{\alpha_k}(\rho;\lambda^{(n)}) =
\lambda_k^*\frac{\alpha_k}{\alpha}(\rho;\lambda^{*})
\frac{\alpha}{\alpha_k}(\rho;\lambda^{(n)}) < \lambda_k^*\frac{\alpha_k}{\alpha}(\rho;\lambda^{(n)})
\frac{\alpha}{\alpha_k}(\rho;\lambda^{(n)}) =  \lambda_k^{*}\,,
\end{equation}
the inequality following for the last time from the inductive hypothesis and from 
properties $p_2$) and $p_3$) of Proposition~3.1.
\end{proof}
According to Propositions 4.1 and 4.2, the sequence converges. Now, let $\tilde \lambda =
\lim_{n\to\infty}\,\lambda^{(n)}$ be the limit of the sequence. Effortlessly,
we prove that $\tilde\lambda$ is a fixed point of $T$. Indeed, 
\begin{equation}
\tilde\lambda_k \,=\, \lim_{n\to\infty}\lambda_k^{(n)} \,=\,
\lim_{n\to\infty}\mu_k\frac{\alpha}{\alpha_k}(\rho;\lambda^{(n-1)}) \,=\,
\mu_k\frac{\alpha}{\alpha_k}\left(\rho;\lim_{n\to\infty}\lambda^{(n-1)}\right) \,=\, T_k(\tilde\lambda)\,.
\end{equation}
Note that passing the limit over $n$ under the integral sign is certainly
allowed for Gaussian integrals. 
\vskip 0.4cm
\begin{prop}[Uniqueness of the fixed point]
Let $\lambda' = T(\lambda')$ and $\lambda'' = T(\lambda'')$ be two
fixed points of $T$, corresponding to the same choice of $v$, $\rho$ and
$\mu\in\cD(\tau_\rho^{-1})$. Then, it must be $\lambda' = \lambda''$. 
\end{prop}
\begin{proof}
According to the hypothesis, $\lambda'$ and $\lambda''$ fulfill the equations
\begin{align}
\lambda'_k & = \mu_k\frac{\alpha}{\alpha_k}(\rho;\lambda') \qquad\hskip 0.5ex \Rightarrow \qquad
\mu_k = \lambda'_k\frac{\alpha_k}{\alpha}(\rho;\lambda')\,,\\[2.0ex]
\lambda''_k & = \mu_k\frac{\alpha}{\alpha_k}(\rho;\lambda'') \qquad \Rightarrow \qquad
\mu_k = \lambda''_k\frac{\alpha_k}{\alpha}(\rho;\lambda'')\,.
\end{align}
Hence, 
\begin{equation}
0 = \mu_k - \mu_k = \lambda'_k\frac{\alpha_k}{\alpha}(\rho;\lambda') -
\lambda_k''\frac{\alpha_k}{\alpha}(\rho;\lambda'') = \sum_{\ell=1}^v \left[\int_0^1\rd t\
J_{k\ell}\left(\rho;\lambda'' + t\left(\lambda'-\lambda''\right)\right)\right](\lambda'_\ell - \lambda''_\ell)\,,
\label{eq:fixedpointid}
\end{equation}
where $J$ denotes the Jacobian matrix of $\tau_\rho$ and is given by
\begin{equation}
J_{k\ell}(\rho;\lambda)
= \partial_k\left(\lambda_\ell\frac{\alpha_\ell}{\alpha}(\rho;\lambda)\right) =
\frac{1}{2}\frac{\lambda_\ell}{\lambda_k}\left(\frac{\alpha_{k\ell}}{\alpha} - \frac{\alpha_k\alpha_\ell}{\alpha^2}\right) = \left[\Lambda^{-1}\,\Omega(\rho;\lambda)\,\Lambda\right]_{k\ell}\,,
\end{equation}
having set $\Omega_{k\ell}\equiv (1/2)(\alpha_{k\ell}/\alpha -
\alpha_k\alpha_\ell/\alpha^2)$. It will be noted that $\Omega =
\{\Omega_{k\ell}\}_{k,\ell=1}^v$ is essentially the covariance matrix of the
square components of $X$ under spherical truncation (we have come across its
matrix elements in eqs.~(\ref{eq:proofp1})--(\ref{eq:covineq})). 
As such, $\Omega$ is symmetric positive definite. Indeed,
\begin{align}
\Omega_{k\ell} & =\,
\frac{1}{2\lambda_k\lambda_\ell}\cov\left(X_k^2,X_\ell^2\,|\,X\in\cB_v(\rho)\right)
\nonumber\\[1.0ex] 
& = \, \frac{1}{2\lambda_k\lambda_\ell}\E\biggl[\biggl(X_k^2-\E\left[X_k^2\,|\,X\in\cB_v(\rho)\right]\biggr)\biggl(X_\ell^2-\E[X_\ell^2\,|\,X\in\cB_v(\rho)]\biggr)\,\biggr|\,X\in\cB_v(\rho)\biggr]\,.
\end{align}
On setting $Z_k = (X_k^2-\E[X_k^2\,|\,X\in\cB_v(\rho)])/\sqrt{2}\lambda_k$, we can represent
$\Omega$ as $\Omega = \E[Z\trans{Z}\,|\,X\in\cB_v(\rho)]$. If $x\in\RR^v$ is not the null vector,
then $\trans{x}\Omega x = \E[\trans{x}Z\trans{Z}x\,|\,X\in\cB_v(\rho)] =
\E[(\trans{x}Z)^2\,|\,X\in\cB_v(\rho)]>0$. Moreover, the eigenvalues of $\Omega$ 
fulfill the secular equation 
\begin{equation}
0=\det(\Omega - \phi\mathds{I}_v) = \det[\Lambda^{-1}(\Omega -
\phi\mathds{I}_v)\Lambda] = \det(\Lambda^{-1}\Omega\Lambda-\phi\mathds{I}_v) = \det(J-\phi\mathds{I}_v)\,,
\end{equation}
whence it follows that $J$ is positive definite as well (though it is not
symmetric). Since the sum of positive definite matrices is positive definite,
we conclude that $\int_0^1\rd t\ J\left(\rho;\lambda'' +
  t\left(\lambda'-\lambda''\right)\right)$ is positive definite too. As such, it
is non--singular. Therefore, from eq.~(\ref{eq:fixedpointid}) we conclude that
$\lambda'=\lambda''$.  
\end{proof}

\section{Numerical computation of Gaussian integrals over $\cB_v(\rho)$}

Let us now see how to compute $\alpha_{k\ell m\ldots}$ with controlled
precision. Most the relevant work has been originally done by
Ruben in ref.~\cite{ruben4}, where the case of $\alpha$ is discussed. 
We extend Ruben's technique to Gaussian integrals containing powers of the integration variable. 
Specifically, it is shown in ref.~\cite{ruben4} that $\alpha(\rho;\lambda)$ can be 
represented as a series of chi--square cumulative distribution functions, 
\begin{equation}
\alpha(\rho;\lambda) = \sum_{m=0}^\infty c_m(s;\lambda)F_{v+2m}(\rho/s)\,.
\label{eq:ruben0}
\end{equation}
The scale factor $s$ has the same physical dimension as $\rho$ and
$\lambda$. It is introduced in order to factorize the
dependence of $\alpha$ upon $\rho$ and $\lambda$ at each order of the expansion. The
series on the {\it r.h.s.} of eq.~(\ref{eq:ruben0}) converges uniformly on
every finite interval of $\rho$. The coefficients $c_m$ are given by
\begin{equation}
c_m(s;\lambda) =
\frac{1}{m!}\frac{s^{v/2+m}}{|\Lambda|^{1/2}}\frac{\Gamma(v/2+m)}{\Gamma(v/2)}\M[(-Q)^m]\,,\qquad
m=0,1,\ldots
\label{eq:cm}
\end{equation}
having defined $Q(x) \equiv \trans{x}[\Lambda^{-1} - s^{-1}\mathds{I}_v]x$ for
$x\in\RR^v$ and $\M$ as the uniform average operator on the $(v-1)$--sphere
$\partial\cB_v(1) \equiv \{u\in\RR^v:\ \trans{u}u=1\}$, {\it viz.}
\begin{equation}
\M[\phi] \,\equiv\, \frac{\Gamma(v/2)}{2\pi^{v/2}}\int_{\partial\cB_v(1)}\rd
u\, \phi(u)\,,\qquad \forall\ \phi\in\cC^{0}\left(\partial\cB_v\left(1\right)\right)\ {\rm a.e.}
\label{eq:unifaverop}
\end{equation}
Unfortunately, eq.~(\ref{eq:cm}) is not particularly convenient for numerical
computations, since $\M[(-Q)^m]$ is only given in integral
form. However, it is also shown in ref.~\cite{ruben4} that the coefficients $c_m$
can be extracted from the Taylor expansion (at $z_0=0$) of the generating function
\begin{equation}
\psi(z) \, = 
\prod_{k=1}^v\left(\frac{s}{\lambda_k}\right)^{1/2}\left[1-\left(1-\frac{s}{\lambda_k}\right)z\right]^{-1/2}
\,, \qquad {\rm i.e.}\quad \psi(z) = \sum_{m=0}^\infty c_m(s;\lambda)z^m\,.
\end{equation}
This series converges uniformly for
$|z|<\min_i|1-s/\lambda_i|^{-1}$. On evaluating the derivatives of $\psi(z)$,
it is then shown that the $c_m$'s fulfill the recursion
\begin{equation}
\left\{\begin{array}{l}
\displaystyle{c_0 \,=\, \prod_{m=1}^v\sqrt\frac{s}{\lambda_m}\,; \qquad c_n =
  \frac{1}{2n}\sum_{r=0}^{n-1}g_{n-r}c_r}\,,\qquad n=1,2,\ldots\,;\\[4.0ex]
\displaystyle{g_n \,\equiv\, \sum_{m=1}^v\left(1 - \frac{s}{\lambda_m}\right)^{n}}\,.
\end{array}\right.
\end{equation}
Finally, the systematic error produced on considering only the lowest $k$
terms of the chi--square series of eq.~(\ref{eq:ruben0}) is estimated by 
\begin{align}
\cR_n(\rho;\lambda) & \,\equiv\,\left|\sum_{m=n}^\infty
  c_m(s;\lambda)F_{v+2m}(\rho/s)\right| \nonumber\\[2.0ex] & \le
c_0(s;\lambda)\frac{\Gamma(v/2+n)}{\Gamma(v/2)}
\frac{\eta^n}{n!}(1-\eta)^{-(v/2+n)}F_{v+2n}[(1-\eta)\rho/s] \,\equiv\,{\frak R}_n\,,
\label{eq:errest}
\end{align}
with $\eta = \max_{i}|1-s/\lambda_i|$. 

Now, as mentioned, it is possible to extend the above expansion to all
Gaussian integrals $\alpha_{k\ell m\ldots}$. Here, we are interested only in
$\alpha_k$ and $\alpha_{jk}$, since these are needed in order to implement the
fixed point iteration and to compute the Jacobian matrix of $\tau_\rho$. The
extension is provided by the following  
\begin{theo}[Ruben's expansions]
The integrals $\alpha_k$ and $\alpha_{jk}$ admit the series representations
\begin{align}
\label{eq:ruben1}\alpha_k(\rho;\lambda) & \,=\, \sum_{m=0}^{\infty}c_{k;m}(s;\lambda)F_{v+2(m+1)}(\rho/s)\,, \\[1.0ex]
\label{eq:ruben2}\alpha_{jk}(\rho;\lambda) & \,=\, \sum_{m=0}^{\infty}c_{jk;m}(s;\lambda)F_{v+2(m+2)}(\rho/s)\,,
\end{align}
with $s$ an arbitrary positive constant. The series coefficients are
given resp. by
\begin{align}
\label{eq:rcoefs1}
c_{k;m}(s;\lambda) & \,=\,
\frac{s}{\lambda_k}\,\frac{v+2m}{m!}\,\frac{s^{v/2+m}}{|\Lambda|^{1/2}}\,\frac{\Gamma(v/2+m)}{\Gamma(v/2)}\,\M\left[(-Q)^m
u_k^2\right]\,,
\end{align}
\begin{align}
\label{eq:rcoefs2}
c_{jk;m}(s;\lambda) & \,=\,
(1+2\delta_{jk})\frac{s}{\lambda_j}\frac{s}{\lambda_k}\,\frac{(v+2m+2)(v+2m)}{m!}\,\frac{s^{v/2+m}}{|\Lambda|^{1/2}}\,\frac{\Gamma(v/2+m)}{\Gamma(v/2)}\,\M\left[(-Q)^m
u_j^2u_k^2\right]\,.
\end{align}
with $\delta_{jk}$ denoting the Kronecker symbol. The series on the r.h.s. of eqs.~(\ref{eq:ruben1})--(\ref{eq:ruben2})
converge uniformly on every finite interval of $\rho$. The functions
\begin{align}
\label{eq:genfunck}
\psi_k(z) & \, = 
\left(\frac{s}{\lambda_{k}}\right)^{3/2}\left[1-\left(1-\frac{s}{\lambda_k}\right)z\right]^{-3/2}\prod_{i\ne
k}\left(\frac{s}{\lambda_i}\right)^{1/2}\left[1-\left(1-\frac{s}{\lambda_i}\right)z\right]^{-1/2}\,,\\[1.0ex]
\label{eq:genfunckk}
\psi_{kk}(z) & \, =  3\left(\frac{s}{\lambda_{k}}\right)^{5/2}\left[1-\left(1-\frac{s}{\lambda_k}\right)z\right]^{-5/2}\prod_{i\ne
k}\left(\frac{s}{\lambda_i}\right)^{1/2}\left[1-\left(1-\frac{s}{\lambda_i}\right)z\right]^{-1/2}\,,
\end{align}
\begin{align}
\label{eq:genfuncjk}
\psi_{jk}(z) & \, = \,
\left(\frac{s}{\lambda_{j}}\frac{s}{\lambda_{k}}\right)^{3/2}\left\{\left[1-\left(1-\frac{s}{\lambda_j}\right)z\right]
\left[1-\left(1-\frac{s}{\lambda_k}\right)z\right]\right\}^{-3/2}\nonumber\\[1.0ex]
&\, \times\prod_{i\ne
j,k}\left(\frac{s}{\lambda_i}\right)^{1/2}\left[1-\left(1-\frac{s}{\lambda_i}\right)z\right]^{-1/2}\
\quad (j\ne k)\,,
\end{align}
are generating functions resp. for the coefficients $c_{k;m}$,
$c_{kk;m}$ and $c_{jk;m}$ ($j\ne k$), i.e. they fulfill
\begin{align}
\psi_k(z) & = \sum_{m=0}^\infty c_{k;m}(s;\lambda)z^m\,,\\[1.0ex]
\psi_{jk}(z) & = \sum_{m=0}^\infty c_{jk;m}(s;\lambda)z^m\,,
\end{align}
for $|z|<\min_i|1-s/\lambda_i|^{-1}$. Finally, the coefficients $c_{k;m}$, $c_{kk;m}$ and $c_{jk;m}$ ($j\ne k$) can
be obtained iteratively from the recursions
\begin{equation}
\left\{\begin{array}{l}
\displaystyle{c_{k;0} \,=\, \left(\frac{s}{\lambda_k}\right)c_0\,; \qquad c_{k;m} =
  \frac{1}{2m}\sum_{r=0}^{m-1}g_{k;m-r}\,c_{k;r}}\,;\\[2.0ex]
\displaystyle{g_{k;m} \,\equiv\, \sum_{i=1}^v e_{k;i}\left(1 - \frac{s}{\lambda_i}\right)^{m}}\,,\qquad m\ge 1\,;
\end{array}\right.
\label{eq:recurs1}
\end{equation}
and
\begin{equation}
\left\{\begin{array}{l}
\displaystyle{c_{jk;0} \,=\, (1+2\delta_{jk})\left(\frac{s}{\lambda_j}\right)\left(\frac{s}{\lambda_k}\right)c_0\,; \qquad c_{jk;m} =
  \frac{1}{2m}\sum_{r=0}^{m-1}g_{jk;m-r}\,c_{jk;r}}\,;\\[2.0ex]
\displaystyle{g_{jk;m} \,\equiv\, \sum_{i=1}^v e_{jk;i}\left(1 - \frac{s}{\lambda_i}\right)^{m}}\,,\qquad m\ge 1\,;
\end{array}\right.
\end{equation}
where the auxiliary coefficients $e_{k;i}$ and $e_{jk;i}$ are defined by
\begin{align}
e_{k;i} & = \left\{\begin{array}{ll}
3 & \text{if \ } i=k\\[1.0ex]
1 & \text{otherwise}
\label{eq:auxcoefs1}
\end{array}\right.\,,\\[2.0ex]
e_{kk;i} & = \left\{\begin{array}{ll}
5 & \text{if \ } i=k\\[1.0ex]
1 & \text{otherwise}
\end{array}\right.,
\quad\qquad
e_{jk;i} = \left\{\begin{array}{ll}
3 & \text{if \ } i= j\ or\ k\\[1.0ex]
1 & \text{otherwise}
\end{array}\right.\ (j\ne k)\,.
\label{eq:auxcoefs2}
\end{align}
\end{theo}
\hfill\qedsymbol

It is not difficult to further generalize this theorem, so as to provide
a chi--square expansion for any Gaussian integral $\alpha_{k\ell m\ldots}$. The
proof follows closely the original one given by Ruben. We reproduce it in
Appendix A for $\alpha_k$, just to highlight the differences
arising when the Gaussian integral contains powers of the integration
variable. 

Analogously to eq.~(\ref{eq:errest}), it is possible to estimate the
systematic error produced when considering only the lowest $k$ terms of the
chi--square series of $\alpha_k$ and $\alpha_{jk}$. Specifically, we find 
\begin{align}
\cR_{k;n} & \,\equiv\,  \left|\sum_{m=n}^\infty
  c_{k;m}(s;\lambda)F_{v+2(m+1)}(\rho/s)\right| \nonumber\\[2.0ex] & \le\,
c_{k;0}\frac{\eta^n}{n!}(1-\eta)^{-(v/2+n+1)}\frac{\Gamma(v/2+n+1)}{\Gamma(v/2)}F_{v+2(n+1)}[(1-\eta)\rho/s]\,\equiv\,{\frak
R}_{k;n}\,,
\end{align}
\begin{align}
\cR_{jk;n} & \,\equiv\, \left|\sum_{m=n}^\infty
  c_{jk;m}(s;\lambda)F_{v+2(m+2)}(\rho/s)\right| \nonumber\\[2.0ex] & \le
c_{jk;0}\frac{\eta^n}{n!}(1-\eta)^{-(v/2+n+2)}\frac{\Gamma(v/2+n+2)}{\Gamma(v/2)}F_{v+2(n+2)}[(1-\eta)\rho/s]\,\equiv\,{\frak
R}_{jk;n}\,.
\end{align}
In order to evaluate all Ruben series with controlled uncertainty, we first
set\footnote{see once more ref.~\cite{ruben4} for an exhaustive discussion on how to
  choose $s$.} $s=2\lambda_1\lambda_v/(\lambda_1+\lambda_v)$, then we choose a 
unique threshold $\varepsilon$ representing the maximum tolerable systematic error,
{\it e.g.} $\varepsilon_{\rm dp} = 1.0\times 10^{-14}$ (roughly corresponding
to double floating--point precision), for all $\alpha$, $\alpha_k$ and
$\alpha_{jk}$, and finally for each $\alpha_X$ we compute the integer
\begin{equation}
\kcut \,\equiv\,\min_{n\ge 1}\left\{n:\ {\frak
    R}_{X;n}<\varepsilon\right\}\,,
\end{equation}
providing the minimum number of chi-square terms, for which the upper bound
${\frak R}_{X;n}$ to the residual sum $\cR_{X;n}$ lies below $\varepsilon$. Of
course, this procedure overshoots the minimum number of terms 
really required for the $\cR$'s to lie below $\varepsilon$, since we actually
operate on the $\frak R$'s instead of the $\cR$'s. Nevertheless, the
computational overhead is acceptable, as it will be shown in next
section. For the sake of completeness, it must be said that typically the
values of $\kcut$ for $\alpha$, $\alpha_k$ and $\alpha_{jk}$ with
the same $\epsilon$ (and $\rho$, $\lambda$) are not much different from each other. 

To conclude, we notice that $\kcut$ depends non--trivially upon
$\lambda$. By contrast, since $F_v(x)$ is monotonic increasing in $x$, we 
clearly see that $\kcut$ is monotonic increasing in $\rho$. Now, should one
evaluate $\alpha$ and the like for a given $\lambda$ at several values of $\rho$, say  
 $\rho_1 \le \rho_2\le \ldots\le \rho_{\rm max}$, it is
advisable to save computing resources and work out Ruben coefficients
just once, up to the order $\kcut$ corresponding to $\rho_{\rm \max}$, since
$\kcut(\rho_1)\le\ldots\le\kcut(\rho_{\rm max})$. We made use of this trick
throughout our numerical experiences, as reported in the sequel.

\section{Numerical analysis of the reconstruction process}

The fixed point eq.~(\ref{eq:fixedpoint}) represents the simplest iterative
scheme that can be used in order to reconstruct the solution $\lambda =
\tau_\rho^{-1}\cdot\mu$. In the literature of numerical methods, this scheme
is known as a non--linear Gauss--Jacobi (GJ) iteration (see {\it e.g.} ref.~\cite{judd}). Accordingly, 
we shall rewrite it as $\lambda^{(n+1)}_{\GJ,k} = T_k(\lambda^{(n)}_{\GJ})$. As we 
have seen, the sequence $\lambda^{(n)}_{\GJ}$ converges with no exception as
$n\to\infty$, provided $\mu\in\cD(\tau_\rho^{-1})$. 
Given $\epsT>0$, the number of steps $\nit$ needed for
an approximate convergence with relative precision $\epsT$, i.e.
\begin{equation} 
\nit \,\equiv\, \min_{n\ge 1}\left\{n:\
  \frac{||\lambda_\GJ^{(n)}-\lambda_\GJ^{(n-1)}||_{\infty}}{||\lambda_\GJ^{(n-1)}||_{\infty}}<\epsT\right\}\,,
\label{eq:nitdef}
\end{equation} 
depends not only upon $\epsT$, but also on $\rho$ and $\mu$ (note that
the stopping rule is well conditioned, since
\nohyphens{$||\lambda^{(n)}||_\infty>~0$} $\forall\,n$ and also $\lim_{n\to
  \infty}||\lambda^{(n)}||_\infty>0$). In order to  
characterize statistically the convergence rate of the reconstruction process, we
must integrate out the fluctuations of $\nit$ due to changes of $\mu$, i.e. we
must average $\nit$ by letting $\mu$ fluctuate across its own probability
space. In this way, we obtain the quantity $\barnit
~\equiv~\E_\mu[\nit|\epsT,\rho]$, which better synthesizes the cost of the reconstruction 
for given $\epsT$ and $\rho$. It should be evident that carrying out this idea analytically is hard,
for on the one hand $\nit$ depends upon $\mu$  non--linearly, and on the other 
$\mu$ has a complicate distribution, as we briefly explain below.  

\subsection{Choice of the eigenvalue ensemble}

Since $\lambda$ is the eigenvalue spectrum of a full
covariance matrix, it is reasonable to assume its distribution to be a Wishart
$\cW_v(p,\Sigma_0)$ for some scale matrix $\Sigma_0$ and for some number of
degrees of freedom $p\ge v$. In the sequel, we shall make the ideal assumption
$\Sigma_0=\idp$, so that the probability measure of $\lambda$ is (see
{\it e.g.} ref.~\cite{anderson}) 
\begin{equation} 
\rd w_v(p;\lambda) \, = \,
p^{(p+v^2-1)/2}\,\frac{\pi^{v^2/2}\prod_{k=1}^v\lambda_k^{(p-v-1)/2}\exp
\left(-\frac{p}{2}\sum_{k=1}^v\lambda_k\right)\prod_{k<j}(\lambda_j-\lambda_k)}
{2^{vp/2}\Gamma_v(p/2)\Gamma_v(v/2)}\,\rd^v\lambda\,.
\label{eq:wishartdistr}
\end{equation}
Under this assumption, the probability measure of $\mu$ is obtained 
by performing the change of variable $\lambda = \tau_\rho^{-1}\cdot\mu$ in
eq.~(\ref{eq:wishartdistr}). Unfortunately, we have no analytic representation
of $\tau_\rho^{-1}$. Thus, we have neither an expression for the distribution of
$\mu$. However, $\mu$ can be extracted numerically as follows:
\begin{itemize}
\item[{\it i})]{generate randomly
$\Sigma\sim\cW_v(p,\idp)$ by means of the Bartlett decomposition
\cite{bartlett};}
\item[{\it ii})]{take the ordered eigenvalue spectrum $\lambda$ of
$\Sigma$;}
\item[{\it iii})]{obtain $\mu$ by applying the truncation operator $\tau_\rho$
    to $\lambda$.} 
\end{itemize}
Note that since $\cW_v(p,\idp)$ is only defined for $p\ge v$, we need to rescale $p$
as $v$ increases. The simplest choice is to keep the ratio $p/v$ fixed. 
The larger this ratio, the closer $\Sigma$ fluctuates around $\mathds{I}_v$
(recall that if $\Sigma\sim\cW_v(p,p^{-1}\cdot\mathds{I}_v)$, then 
$\E[\Sigma_{ij}]=\delta_{ij}$ and $\var(\Sigma_{ij})=p^{-1}[1+\delta_{ij}]$). 
In view of this, large values of $p/v$ are to be avoided, since they reduce 
the probability of testing the fixed point iteration on eigenvalue spectra 
characterized by large condition numbers $n_{\rm cond} \equiv \lambda_v/\lambda_1$. 
For this reason, we have set $p=2v$ in our numerical study.

Having specified an ensemble of matrices from which to extract the eigenvalue
spectra, we are now ready to perform numerical simulations. To begin with, we
report in Fig.~\ref{fig:distro1} the marginal probability density function of
the ordered eigenvalues $\{\lambda_k\}_{k=1}^v$ and their truncated counterparts
$\{\mu_k\}_{k=1}^v$ for the Wishart ensemble
$\cW_{10}(20,20^{-1}\cdot\mathds{I}_{10})$ at $\rho=1$, as obtained
numerically from a rather large sample of matrices ($\simeq 10^6$ units). It will be
noted that {\it i}) the effect of the truncation is severe on the largest
eigenvalues, as a consequence of the analytic bounds of Corollary~2.1 and
Proposition~2.2; {\it ii}) while the skewness of the lowest truncated
eigenvalues is negative, it becomes positive for the largest ones. This is due
to a change of relative effectiveness of eq.~(\ref{eq:mubounds2}) {\it i})
with respect to eq.~(\ref{eq:mubounds2})~{\it ii}). 

\subsection{Choice of the simulation parameters}

In order to explore the dependence of $\barnit$ upon $\rho$, we need to choose
one or more simulation points for the latter. Ideally, it is possible to identify
three different regimes in our problem: $\rho\lesssim\lambda_1$ (strong truncation regime),
$\lambda_1\lesssim\rho\lesssim\lambda_v$ (crossover) and
$\rho\gtrsim\lambda_v$ (weak truncation regime). We cover all of them with the
following set of points:
\begin{equation}
\rho\,\in\,\left\{\Mo{\lambda_1},\ \ldots\ ,\,\Mo{\lambda_v}\right\}\,\cup\,\left\{\frac{1}{2}\Mo{\lambda_1},\,2\,\Mo{\lambda_v}\right\}\,,
\label{eq:simpoints}
\end{equation}
\begin{table}[!t]
  \small
  \begin{center}
    \begin{tabular}{c|cccccccc}
      \hline\hline\\[-2.0ex]
      & $v=3$ & $v=4$ & $v=5$ & $v=6$ & $v=7$ & $v=8$ & $v=9$ & $v=10$   \\ \\[-2.0ex]
      \hline\\[-1.0ex]
     $\hatMo{\lambda_1}$    & 0.1568 & 0.1487 & 0.1435 & 0.1383 & 0.1344 & 0.1310 & 0.1269 & 0.1258 \\[2.0ex]
     $\hatMo{\lambda_2}$    & 0.6724 & 0.4921 & 0.4017 & 0.3424 & 0.3039 & 0.2745 & 0.2554 & 0.2399 \\[2.0ex]
     $\hatMo{\lambda_3}$    & 1.6671 & 1.0112 & 0.7528 & 0.6071 & 0.5138 & 0.4543 & 0.4048 & 0.3693 \\[2.0ex]
     $\hatMo{\lambda_4}$    & --    &  1.8507 & 1.2401 & 0.9621 & 0.7854 & 0.6684 & 0.5858 & 0.5288 \\[2.0ex]
     $\hatMo{\lambda_5}$    & --    & --      & 2.0150 & 1.4434 & 1.1269 & 0.9263 & 0.7956 & 0.7032 \\[2.0ex]
     $\hatMo{\lambda_6}$    & --    & --    & --       & 2.1356 & 1.5789 & 1.2559 & 1.0527 & 0.9111 \\[2.0ex]
     $\hatMo{\lambda_7}$    & --    & --    & --       & --     & 2.2190 & 1.6764 & 1.3673 & 1.1603 \\[2.0ex]
     $\hatMo{\lambda_8}$    & --    & --    & --       & --     & --     & 2.2763 & 1.7687 & 1.4624 \\[2.0ex]
     $\hatMo{\lambda_9}$    & --    & --    & --       & --     & --     & --     & 2.3210 & 1.8473 \\[2.0ex]
     $\hatMo{\lambda_{10}}$ & --    & --    & --       & --     & --     & --     & --     & 2.3775 \\[2.0ex]

      \hline\hline
    \end{tabular}
    \vskip 0.4cm
    \caption{Numerical estimates of the mode of the ordered eigenvalues
      $\{\lambda_1,\dots,\lambda_v\}$ of
      $\Sigma\sim\cW_v\left(2v,(2v)^{-1}\cdot\mathds{I}_v\right)$ with
      $v=3,\dots,10$. The estimates have been obtained from Grenander's
      mode estimator \cite{grenander}.\label{tab:modes}}  
  \end{center}
    \vskip -0.5cm
\end{table}

\vskip -0.4cm
\noindent where $\Mo{\cdot}$ stands for the mode. In principle, it is possible to
determine $\Mo{\lambda_k}$ with high accuracy by using analytic
representations of the marginal probability densities of the ordered
eigenvalues~\cite{zanella}. In practice, the latter become computationally
demanding at increasingly large values of $v$: for instance, the determination of the
probability density of $\lambda_2$ requires $(v!)^2$ sums, which is unfeasible
even at $v\sim 10$. Moreover, to our aims it is sufficient to  choose
approximate values, provided these lie not far from the exact ones.
Accordingly, we have determined the eigenvalue modes numerically from samples
of $N\simeq 10^6$  Wishart matrices. Our estimates are reported in
Table~\ref{tab:modes} for $v=3,\ldots,10$.  
They have been obtained from Grenander's estimator~\cite{grenander}, 
\begin{equation}
\Mo{\lambda_k}_{rs} = \frac{1}{2}\frac{\sum_{i=1}^{N-r}\left(\lambda_k^{(i)}+\lambda_k^{(i+r)}\right)\left(\lambda_k^{(i)}-\lambda_k^{(i+r)}\right)^{-s}}{\sum_{i=1}^{N-r}\left(\lambda_k^{(i)}-\lambda_k^{(i+r)}\right)^{-s}}\,,
\end{equation}
with properly chosen parameters $r$, $s$. 

\begin{center}
\begin{figure}
\begin{center}
\includegraphics[width=0.9\textwidth]{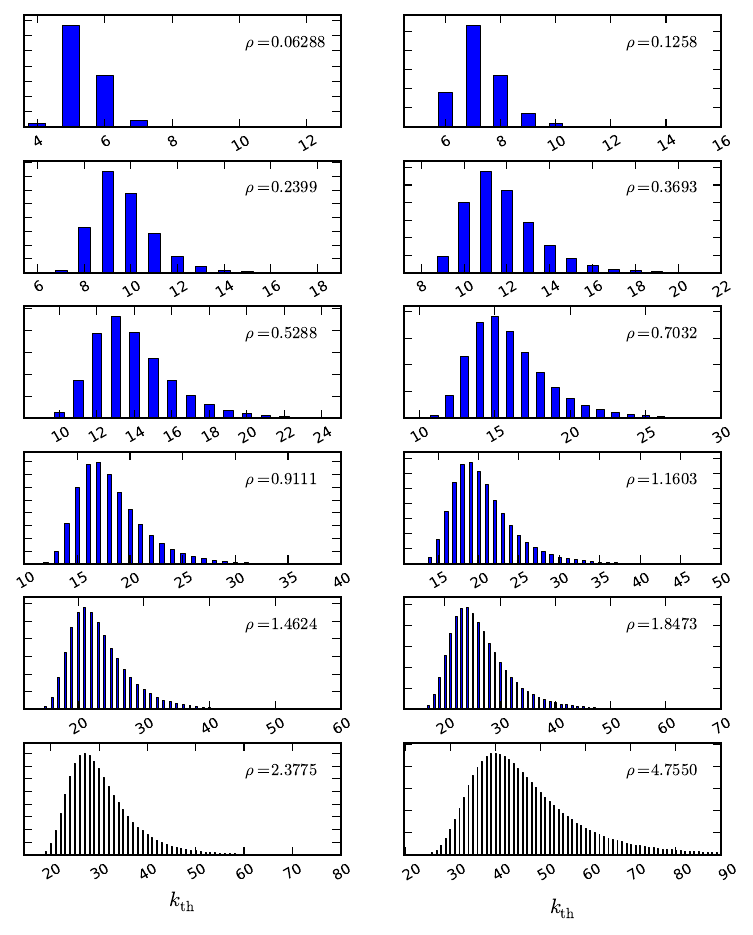}
\vskip -0.2cm
\caption{\small Monte Carlo simulation of the probability mass function of
  the parameter $\kcut$ for the Gaussian probability content
  $\alpha$. The histograms refer to the eigenvalue ensemble
  $\lambda$ of $\Sigma\sim\cW_{10}(20,20^{-1}\cdot\mathds{I}_{10})$, with $\rho$ chosen 
  as in Table~\ref{tab:modes} and $\varepsilon = 1.0\times 10^{-14}$.\label{fig:distro2}} 
\end{center}
\end{figure}
\end{center}

We are now in the position to investigate numerically how many terms in
Ruben's expansions must be considered as $\varepsilon$ is set to
$\varepsilon_{\rm dp}=1.0\times 10^{-14}$, for our choice of the eigenvalue
ensemble $\lambda\sim\cW_v\left(2v,(2v)^{-1}\cdot\mathds{I}_v\right)$ and with $\rho$ set
as in Table~\ref{tab:modes}. As an example, we report in
Fig.~\ref{fig:distro2} the discrete distributions of $\kcut$ for
the basic Gaussian integral $\alpha$ at $v=10$, the largest dimension we have
simulated. As expected, we observe an increase of $\kcut$ with
$\rho$. Nevertheless, we see that the number of Ruben's components
to be taken into account for a double precision result keeps altogether modest
even in the weak truncation regime, which proves the practical usefulness of the
chi--square expansions.  

\subsection{Fixed point iteration at work}

The GJ iteration is too slow to be of practical interest. For instance, 
at $v=10$, $\rho\simeq\Mo{\lambda_1}$ and $\epsT=1.0\times 10^{-7}$ 
(corresponding to a reconstruction of $\lambda$ with single floating--point
precision) it is rather easy to extract realizations of $\mu$ which require $\nit\simeq 15,000$
to converge. An improvement of the GJ scheme is achieved via over--relaxation (GJOR), i.e. 
\begin{equation}
\left\{\begin{array}{l}
\lambda_{\GJOR,k}^{(0)} \, = \, \mu_k\,, \\[2.0ex]
\lambda_{\GJOR,k}^{(n+1)}\, =\, \lambda_{\GJOR,k}^{(n)} + \omega\left[T_k(\lambda_{\GJOR}^{(n)}) - \lambda_{\GJOR,k}^{(n)}\right]\,,\qquad k=1,\ldots,v
\end{array}\right.\,.
\label{eq:GJOR}
\end{equation}
\begin{center}
\begin{figure}
\begin{center}
\includegraphics[width=0.9\textwidth]{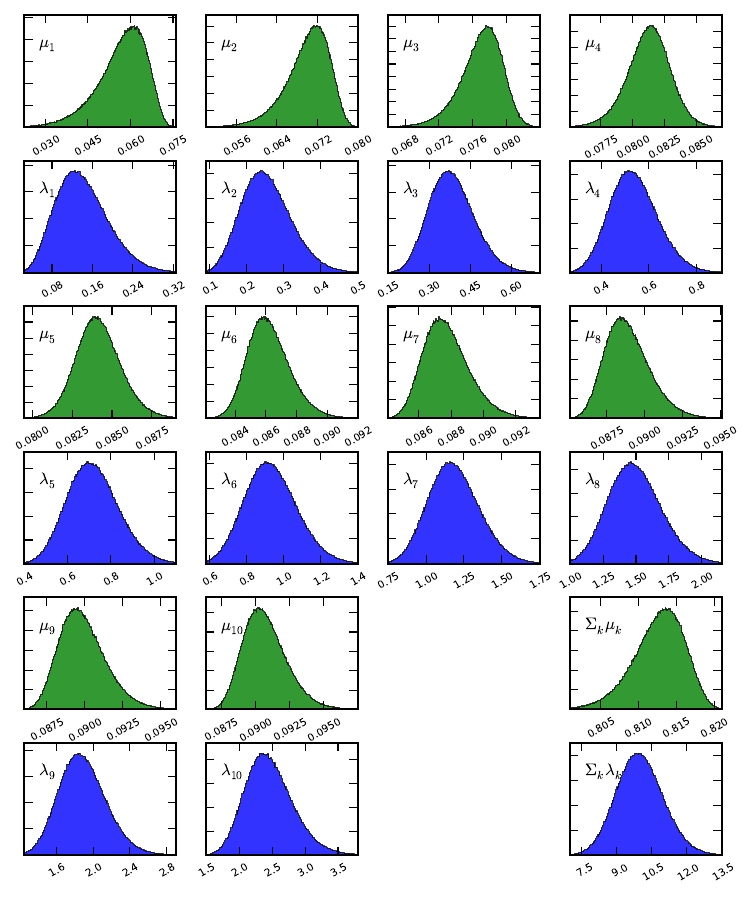}
\vskip -0.2cm
\caption{\small Monte Carlo simulation of the probability density function
  of the ordered eigenvalues $\lambda_k$ (even rows) and their truncated
  counterparts $\mu_k$ at $\rho=1$ (odd rows) for the Wishart ensemble
  $\cW_{10}(20,20^{-1}\cdot\mathds{I}_{10})$. The last two plots (bottom
  right) display the distribution of the sum of
  eigenvalues.\label{fig:distro1}}    
\end{center}
\end{figure}
\end{center}
\vskip -0.3cm
\noindent Evidently, at $\omega=1$ the GJOR scheme coincides with the standard GJ
one.  The optimal value $\omega_{\rm opt}$ of the relaxation factor $\omega$
is not obvious even in the linear Jacobi scheme, where $\omega_{\rm opt}$ depends
upon the properties of the coefficient matrix of the system. For instance, if
the latter is symmetric positive definite, it is demonstrated that the best
choice is provided by $\omega_{\rm opt} \equiv~2(1+\sqrt{1-\sigma^2})^{-1}$,
being $\sigma$ the spectral radius of the Jacobi iteration
matrix~\cite{young}. In our numerical tests with the GJOR scheme, we found
empirically that the optimal value of $\omega$ at $\rho\lesssim\lambda_v$ is
close to the linear prediction, provided $\sigma$ is replaced by
$||J||_\infty$, being $J$ defined as in sect.~3 (note that $||J||_\infty<1$). By
contrast, the iteration diverges after few steps with increasing probability
as $\rho/\lambda_v\to \infty$ if $\omega$ is kept fixed at $\omega =
\omega_{\rm opt}$; in order to restore the convergence, $\omega$ must be
lowered towards $\omega=1$ as such limit is taken.   

\begin{center}
  \begin{figure}[!t]
    \begin{center}
      \begin{tabular}{p{0.45\textwidth}p{0.45\textwidth}}
        \includegraphics[width=0.45\textwidth]{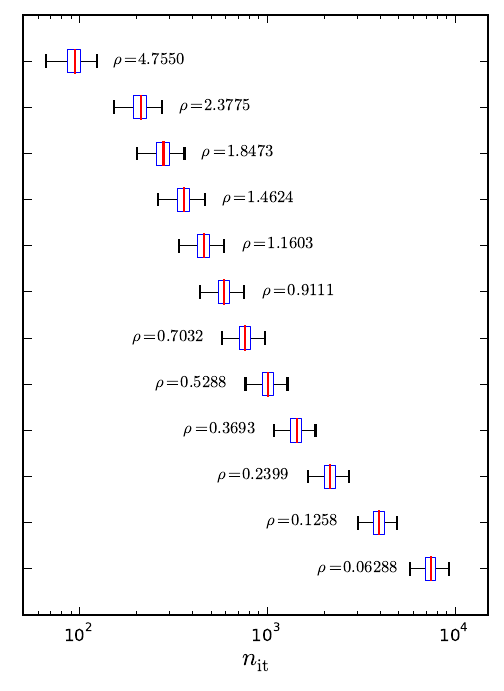}
        \begin{center}
        \end{center}
        &
        \vskip -9.4cm
        \begin{tabular}{c|cc}
          \hline\hline\\[-1.0ex]
          $v$ & \phantom{ide}$\hat a$\,({\it jk.err.})\phantom{ide} &
          \phantom{ide}$\hat b$\,({\it jk.err.})\phantom{ide} \\ \\[-2.0ex]
          \hline\\[-1.0ex]
          3  & 4.93(3) & 0.940(1) \\[2.0ex]
          4  & 5.25(2) & 0.940(1) \\[2.0ex]
          5  & 5.44(2) & 0.948(1) \\[2.0ex]
          6  & 5.65(1) & 0.953(1) \\[2.0ex]
          7  & 5.85(1) & 0.948(1) \\[2.0ex]
          8  & 6.02(1) & 0.951(1) \\[2.0ex]
          9  & 6.18(1) & 0.946(1) \\[2.0ex]
          10 & 6.31(1) & 0.948(1) \\[2.0ex]
          \hline\hline
        \end{tabular}
        \begin{center}
        \end{center}
      \end{tabular}
      \vskip -1.3cm
      \caption{\small {\bf (left)} Box--plot of $\nit$ in the GJOR scheme at
  $v=10$, with $\epsT=1.0\times 10^{-7}$ and $\rho$ chosen as in
  Table~\ref{tab:modes}. The distributions have been reconstructed from a
  sample of $N\simeq 10^3$ eigenvalue spectra extracted from
  $\cW_{10}(20,20^{-1}\cdot\mathds{I}_{10})$. The whiskers extend to the most
  extreme data point within $(3/2)(75\%-25\%)$ data range. {\bf (right)} Numerical
  estimates of the scaling parameters $a$ and $b$ of the GJOR scheme, as obtained
  from jackknife fits to eq.~(\ref{eq:scalinglaw}) of data points with
  $\rho\lesssim 1$ and $\epsT=1.0\times 10^{-7}$. We quote in parentheses the
  jackknife error.\label{fig:sorgjplot}}    
    \end{center}
  \end{figure}
\end{center}

\begin{center}
\begin{figure}
\begin{center}
\includegraphics[width=0.9\textwidth]{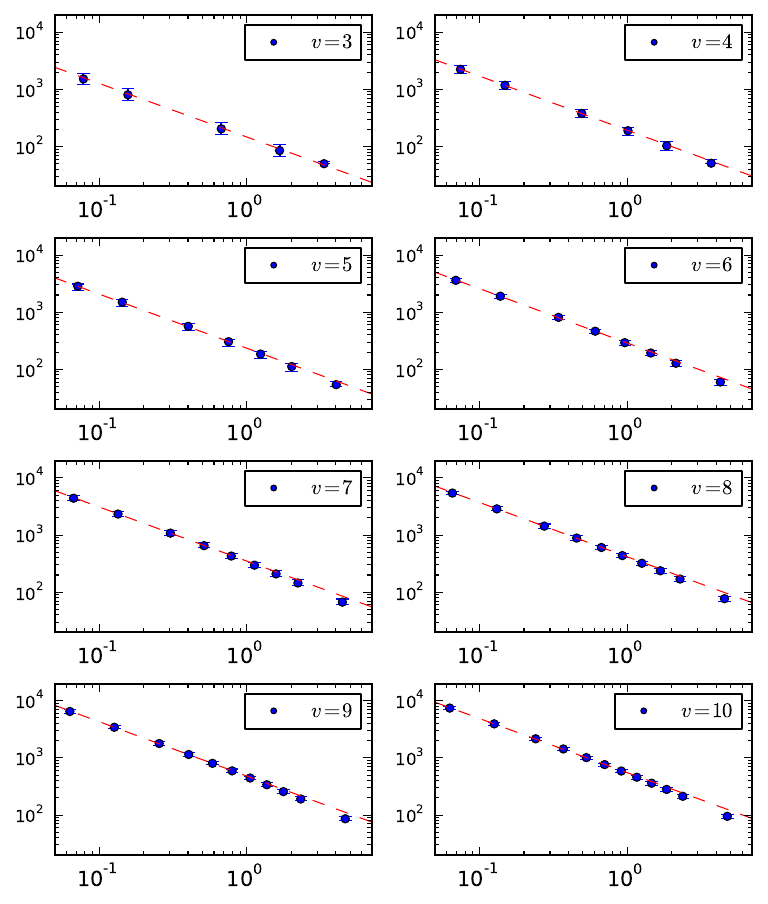}
\vskip -0.2cm
\caption{\small Log--log plots of $\barnit$ vs. $\rho$ in the GJOR scheme at
  $\epsT=1.0\times 10^{-7}$. The parameter $\rho$ has been chosen as in
  Table~\ref{tab:modes}. The (red) dashed line in each plot represents our best
  jackknife linear fit to eq.~(\ref{eq:scalinglaw}) of data points with 
  $\rho\lesssim 1$.\label{fig:nitvsrho}}     
\end{center}
\end{figure}
\end{center}

To give an idea of the convergence rate of the GJOR scheme, we show in
Fig.~\ref{fig:sorgjplot} (left) a joint box--plot of the distributions of
$\nit$ at $v=10$ and $\epsT=1.0\times 10^{-7}$. From the plot we observe
that the distribution of $\nit$ shifts rightwards as $\rho$ decreases:
clearly, the reconstruction is faster if $\rho$ is in the weak
truncation regime (where $\mu$ is closer to $\lambda$), whereas it takes more
iterations in the strong truncation regime. The dependence of $\barnit$ upon
$\rho$, systematically displayed in Fig.~\ref{fig:nitvsrho}, is compatible
with a scaling law   
\begin{equation}
\log \barnit(\rho,v,\epsT) = a(v,\epsT) - b(v,\epsT)\log\rho\,,
\label{eq:scalinglaw}
\end{equation}
\vskip -0.6cm
apart from small corrections occurring at large
$\rho$. Eq.~(\ref{eq:scalinglaw}) tells us that $\barnit$ increases
polynomially in $1/\rho$ at fixed $v$. In order to estimate the parameters $a$
and $b$ in the strong truncation regime (where the algorithm becomes
challenging),  
we performed jackknife fits to eq.~(\ref{eq:scalinglaw}) of data points with $\rho\lesssim 1$. 
Results are collected in Fig.~\ref{fig:sorgjplot} (right), showing that $b$ is
roughly constant, while $a$ increases almost linearly in $v$. Thus, while the
cost of the eigenvalue reconstruction is only polynomial in $1/\rho$ at fixed
$v$, it is exponential in $v$ at fixed $\rho$. The scaling law of the GJOR
scheme is therefore better represented by $\barnit = C\re^{\kappa v}/\rho^b$,
with $C$ being a normalization constant independent of $\rho$ and $v$, and
$\kappa$ representing approximately the slope of $a$ as a function of
$v$. Although the GJOR scheme improves the GJ one, the iteration reveals to be
still inefficient in a parameter subspace, which is critical for the
applications.  

\subsection{Boosting the GJOR scheme}

A further improvement can be obtained by letting $\omega$ depend on the
eigenvalue index in the  GJOR scheme. Let us discuss how to work out such an
adjustment. On commenting Fig.~\ref{fig:distro1}, we
have already noticed that the largest eigenvalues are affected by the
truncation to a larger extent than the smallest ones. Therefore, they must
perform a longer run through the fixed point iteration, in 
order to converge to the untruncated values. This is a possible qualitative
explanation for the slowing down of the algorithm as $\rho\to 0$. In view of
it, we expect to observe some acceleration of the convergence rate, if
$\omega$ is replaced, for instance, by
\begin{equation}
\omega \,\rightarrow\,\omega_k\, \equiv \,(1 + \beta\cdot k)\,\omega_{\rm opt}\,,
\qquad \beta\ge 0\,,\qquad k=1,\ldots,v\,.
\label{eq:boostproposal}
\end{equation}
The choice $\beta=0$ corresponds obviously to the standard GJOR scheme. Any
other choice yields $\omega_k>\omega_{\rm opt}$. Therefore, the new scheme is
also expected to display a higher rate of failures than the
GJOR one at $\rho\gg\lambda_v$, for the reason explained in  sect.~5.3. The
component--wise over--relaxation proposed in eq.~(\ref{eq:boostproposal}) is only meant
to enhance the convergence speed in the strong truncation regime and in the
crossover, where the improvement is actually needed.  

In order to confirm this picture, we have explored systematically the effect
of $\beta$ on $\barnit$ by simulating the re\-con\-struction process at
$v=3,\ldots,10$, with $\beta$ varying from 0 to 2 in steps of $1/5$.
First of all, we have observed that the rate of failures at large
$\rho$ is fairly reduced if the first $30\div 50$ iterations are run with
$\omega_k = \omega_{\rm opt}$, and only afterwards $\beta$ is switched on.
Having minimized the failures, we have
checked that for each value of $\beta$, the scaling law assumed in
eq.~(\ref{eq:scalinglaw}) is effectively fulfilled. Then, we have computed 
jackknife estimates of the scaling parameters $a$ and $b$. These are plotted in
Fig.~\ref{fig:scalingpars} as functions of
$v$. Each trajectory (represented by a dashed curve) corresponds to a given
value of $\beta$. Those with darker markers refer to smaller values of $\beta$
and the other way round. From the plots we notice that  
\begin{itemize}
\item[{\it i})]{all the trajectories with $\beta>0$
    lie below the one with $\beta=0$;}\\[-2.0ex]
\item[{\it ii})]{the trajectories of $a$ display a clear increasing
    trend with $v$, yet their slope lessens as $\beta$ increases. By contrast,
    the trajectories of $b$ develop a mild increasing trend with $v$ as
    $\beta$ increases, though this is not strictly monotonic;}\\[-2.0ex]
\item[{\it iii})]{the trajectories of both $a$ and $b$ seem to converge to a limit
    trajectory as $\beta$ increases; we observe a saturation
    phenomenon, which thwarts the benefit of increasing $\beta$ beyond a
    certain threshold close to $\beta_{\rm max}\simeq 2$.}
\end{itemize}
We add that pushing $\beta$ beyond $\beta_{\rm max}$ is counterproductive, as
the rate of failures becomes increasingly relevant in the crossover and
eventually also in the strong truncation regime. By contrast, if
$\beta\lesssim\beta_{\rm max}$ the rate of failures keeps very low for
essentially all simulated values of $\rho$. 
\begin{center}
  \begin{figure}[!t]
    \begin{center}
      \begin{tabular}{p{0.50\textwidth}p{0.50\textwidth}}
        \includegraphics[width=0.45\textwidth]{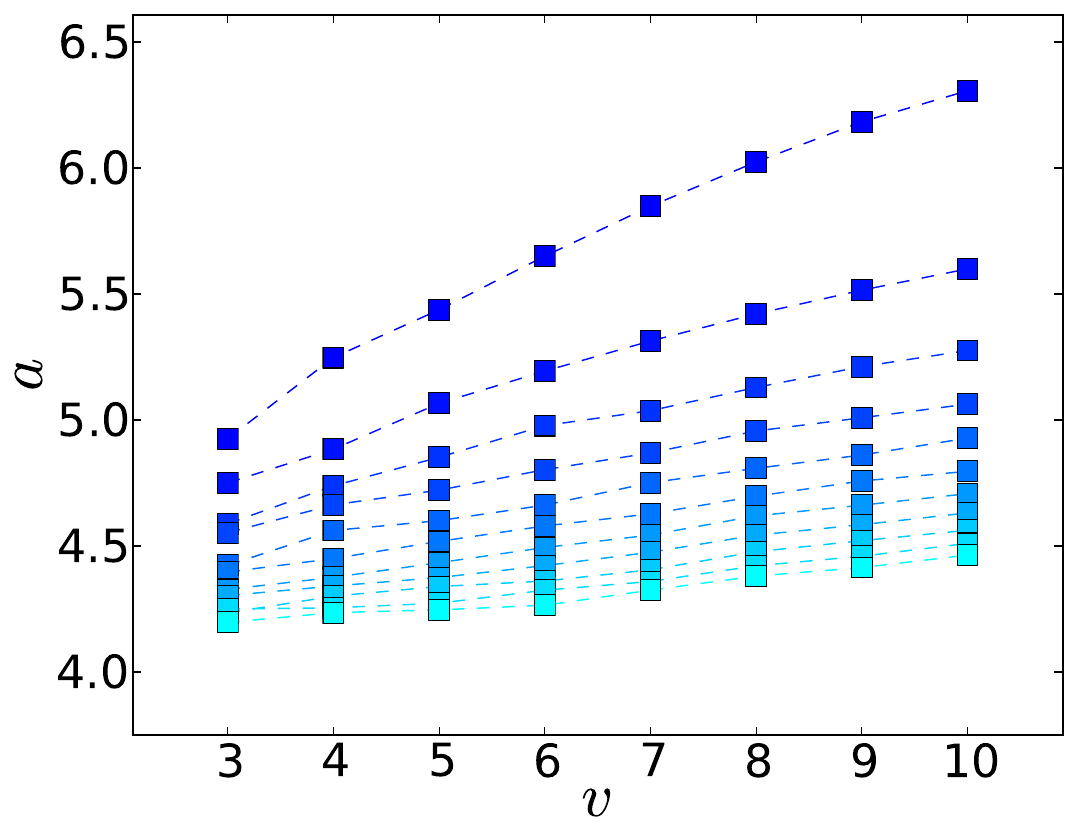}
        \hskip 0.0cm
        &
        \hskip 0.0cm
        \includegraphics[width=0.45\textwidth]{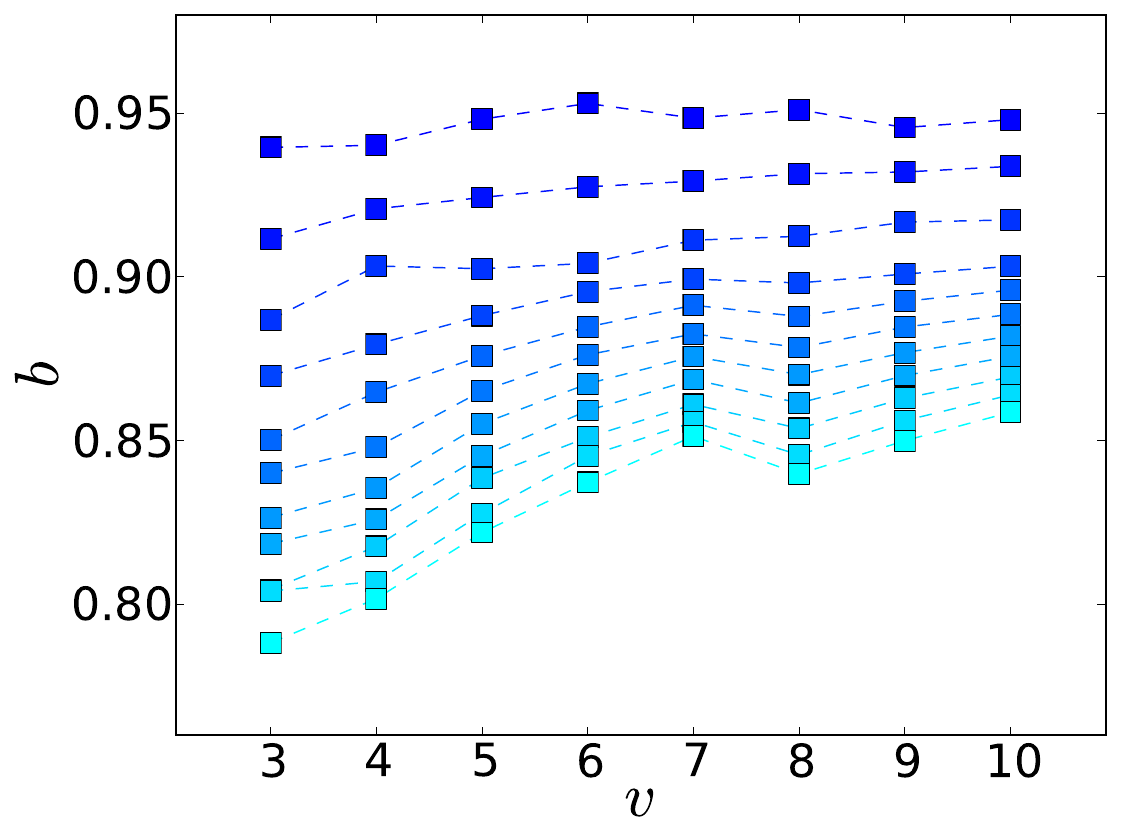}
      \end{tabular}
      \caption{\small Scaling parameters $a$ and $b$ of the
        modified GJOR scheme as functions of $v$ at $\epsT=1.0\times 10^{-7}$, with $\beta$
        varying in the range $0\div 2$ in steps of $1/5$. Each trajectory
        (represented by a dashed curve) refers to a different value of
        $\beta$. Those with darker markers correspond to smaller values of
        $\beta$ and the other way round.   \label{fig:scalingpars}}  
    \end{center}
  \end{figure}
  \vskip 1.5cm
  \begin{figure}[!t]
  \vskip 0.0cm
    \begin{center}
        \includegraphics[width=0.75\textwidth]{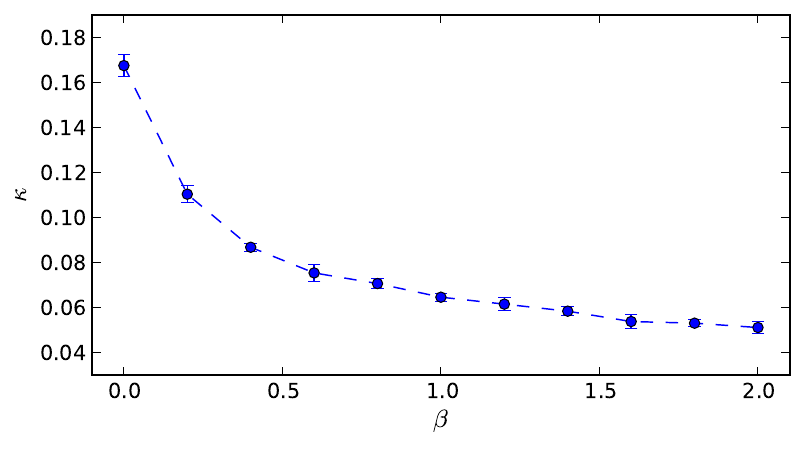}
      \caption{\small The parameter $\kappa$ as a function of $\beta$. Estimates of $\kappa$
        are obtained from least--squares fits of data to a linear model $a =
        a_0 + \kappa\cdot v$.\label{fig:kappafig}}    
    \end{center}
    \vskip -0.3cm
  \end{figure}
\end{center}
\vskip -1.5cm
\noindent Our numerical results signal a strong reduction of the slowing down
of the convergence rate. Indeed, {\it i}) means qualitatively that $C$ and $b$ are
reduced as $\beta$ increases. {\it ii)} means that $\kappa$ 
is reduced as $\beta$ increases (this is the most important effect, as
$\kappa$ is mainly responsible for the exponential slowing down with $v$). The
appearance of a slope in the trajectories of $b$ as $\beta$ increases
indicates that a mild exponential slowing down is also developed at denominator of
the scaling law $\barnit = C\re^{\kappa v}/\rho^b$, but
the value of $b$ is anyway smaller than at $\beta=0$. Finally, {\it iii})
means that choosing $\beta>\beta_{\rm max}$ has a minor impact on the
performance of the algorithm.  In
Fig.~\ref{fig:kappafig}, we report a plot of the parameter $\kappa$ (obtained from
least--squares fits of data to a linear model $a=a_0+\kappa\cdot v$) as a
function of $\beta$. We see that $\kappa(\beta=0)/\kappa(\beta=2)\simeq
4$. This quantifies the maximum exponential speedup of the convergence rate,
which can be achieved by our proposal. When $\beta$ is close to $\beta_{\rm
  max}$, $\barnit$ amounts to few hundreds at $v=10$ and $\rho \simeq
\lambda_1/2$.

\section{On the ill--posedness of the reconstruction in sample space}

So far we have discussed the covariance reconstruction under the
assumption that $\mu = \tau_\rho\cdot\lambda$ represents the exact
truncated counterpart of some $\lambda\in\RR^v$ and we have looked
at the algorithmic properties of the iteration schemes which
operatively define $\tau^{-1}_\rho$. Such analysis is essential in 
order to characterize $\tau^{-1}_\rho$ mathematically, yet it is not
sufficient in real situations, specifically when $\mu$ is perturbed by 
statistical noise. 

In this section, we examine the difficulties arising when 
performing the covariance reconstruction in sample space. We first recall that according 
to Hadamard  \cite{hadamard}, a mathematical problem is well--posed provided 
the following conditions are fulfilled:
\begin{itemize}
  \setlength{\itemindent}{1.0cm}
\item[$H_1$:]{there exists always a solution to the problem;}
\item[\quad $H_2$:]{the solution is unique;}
\item[\quad $H_3$:]{the solution depends smoothly on the input data.}
\end{itemize}
\vskip -0.2cm
Inverse problems are often characterized by violation of one or more of them, see 
for instance ref.~\cite{cavalier}. In such cases, the standard practice consists in 
regularizing the inverse operator, \ie in replacing it by a stable approximation. 
With regard to our problem, the reader will recognize that $H_1$ is violated (and the problem 
becomes ill--posed) as soon as the space of the input data is allowed to be a superset of 
$\cD(\tau_\rho^{-1})$: once clarified how $\mu$ is concretely estimated in the applications 
(sects. 7.1 and 7.2), we propose a perturbative regularization of $\tau_\rho^{-1}$, which 
improves effectively the fulfillment of $H_1$ (sect.~7.3). By contrast, Proposition~4.3 
guarantees that whenever a solution exists, it is also unique, thus $H_2$ is never of concern. 
Finally, the fulfillment of $H_3$ depends on how the statistical noise on $\mu$ is non--linearly 
inflated by the action of $\tau_\rho^{-1}$. For the sake of conciseness, in the present paper 
we just sketch the main ideas underlying the perturbative regularization of $\tau_\rho^{-1}$, 
whereas a technical implementation of it and a discussion of $H_3$ are deferred to a separate 
paper \cite{palopert}.

\subsection{Definition of the sample truncated covariance matrix}

The examples of sect.~2 assume that {\it i})
spherical truncations are operated on a representative sample
${\cal P}_N = \{x^{(k)}\}_{k=1}^N$ of $X\sim{\cal N}_v(0,\Sigma)$
with finite size $N$, {\it ii}) $\rho$ is known exactly and {\it
  iii}) the input budget for the covariance reconstruction is
given by the subset  
\begin{equation}
  {\cal Q}_M = \{ x \in {\cal P}_N:\ ||x||^2<\rho\}\,,\qquad \text{with}\qquad |{\cal
    Q}_M| = M\le N\,.
\end{equation} 
As usual in the analysis of stochastic variables in sample space, we
assume that the observations $x^{(k)}$ are realizations of i.i.d. 
stochastic variables $X^{(k)}\sim{\cal N}_v(0,\Sigma)$,  
$k=1,\ldots,N$. Thus, $M$ is itself a stochastic variable in  
sample space, where it reads
\begin{equation}
M = \sum_{k=1}^{N}\dI_{\cB_v(\rho)}(X^{(k)}) \equiv \sum_{k=1}^N \dI_{k}\,,
\end{equation}
with $\dI_{\cB_v(\rho)}(\cdot)$ denoting the characteristic function of
$\cB_v(\rho)$ and $\dI_k \equiv \dI_{\cB_v(\rho)}(X^{(k)})$ being just a 
shortcut for its extended counterpart. It is easily recognized that $M\sim B(N,\alpha)$ is a binomial variate.
If we indeed denote by $\frak E$ the sample expectation operator (\ie the integral with respect to the product 
measure of the joint variables $\{X^{(k)}\}_{k=1}^N$), then a standard calculation yields
\begin{equation}
\frak{E}[M] =
{\frak E}\left[\sum_{k=1}^{N}\dI_k\right] =
\sum_{k=1}^{N}{\frak E}\left[\dI_k\right] = 
\sum_{k=1}^{N}\alpha = \alpha N\,,
\end{equation}
and
\begin{align}
& {\frak{var}}[M] = {\frak E}[M^2] - {\frak E}[M]^2 \nonumber
= \sum_{k=1}^N {\frak E}\left[\dI_k^2\right] + \sum_{k,s:\ k \ne s}^{1\ldots N}
{\frak E}\left[\dI_k\dI_s\right]  -
\alpha^2N^2 \nonumber\\[0.0ex]
& = \sum_{k=1}^N {\frak E}\left[\dI_k\right] +
\sum_{k,s:\ k\ne s}^{1\ldots N}
{\frak E}\left[\dI_k\right]{\frak E}\left[\dI_s\right]  -
\alpha^2N^2  = \alpha N + \alpha^2N(N-1) -
\alpha^2N^2 = \alpha(1-\alpha)N\,.
\end{align}
Hence, we see that the relative dispersion of $M$ is
$\text{O}(N^{-1/2})$. 
Now, the simplest way to measure $\Sigma$ and ${\frak S}_\cB$
respectively from the sets ${\cal P}_N$ and ${\cal Q}_M$ is
via the classical estimators
\begin{alignat}{3}
\label{eq:samplecovone}
\hat{\Sigma}_{ij} & = \frac{1}{N-1}\sum_{x\in {\cal P}_N} (x-\bar x)_i\,
(x-\bar x)_j\,, & \qquad & \bar x_i = \frac{1}{N}\sum_{x\in{\cal
      P}_N} x_i\,, \\[2.0ex]
(\hat{\frak S}_\cB)_{ij} & = \frac{1}{M-1}\sum_{x\in {\cal Q}_M} (x-\tilde x)_i\,
(x-\tilde x)_j\,, & \qquad & \tilde x_i = \frac{1}{M}\sum_{x\in{\cal Q}_M} x_i\,.
\label{eq:samplecovtwo}
\end{alignat}
We define the sample estimates $\hat\lambda$ and
$\hat\mu$ respectively of $\lambda$ and $\mu$ as the eigenvalue
spectra of $\hat{\Sigma}$ and $\hat{\frak S}_\cB$. By symmetry
arguments we see that $\tilde x_i$ is unbiased. Indeed, it holds  
\begin{equation}
\frak{E}[\tilde x_i] =
\sum_{k=1}^N\frak{E}\left[X_i^{(k)}\dfrac{\dI_k}{\sum_{s=1}^N
  \dI_s}\right]\,.
\label{eq:smeanbias}
\end{equation}
The {\it r.h.s.} of eq.~(\ref{eq:smeanbias})
makes only sense if we conventionally define the integrand to
be zero in the integration subdomain $\{X^{(k)}\notin\cB_v(\rho),\ \forall 
k\}$, or equivalently if we interpret $\frak{E}[\tilde x_i]$ as the conditional one 
$\frak{E}[\tilde x_i\,|\,M>0]$ (the event $M>0$ occurs a.s. only as $N\to\infty$). Since the sample measure is even under $X^{(k)}\to-X^{(k)}$
while the integrand is odd, we immediately conclude that $\frak{bias}[\tilde x_i] =
0$.

\subsection{Bias of the sample truncated covariance matrix}

 The situation gets somewhat less trivial with $\hat {\frak
S}_\cB$: the normalization factor $(M-1)^{-1}$, which has been chosen
in analogy with eq.~(\ref{eq:samplecovone}), is not sufficient to
remove completely the bias of $\hat {\frak S}_\cB$ at finite $N$, though we aim at showing here 
that the residual bias is exponentially small and asymptotically vanishing. In order to see this, 
we observe
\begin{align}
& \frak{E}\left[(\hat{\frak{S}}_\cB)_{ij}\right] = \sum_{k=1}^N
\frak{E}\left[\left(X^{(k)}_i-\tilde X_i\right) \left(X^{(k)}_j-\tilde X_j\right)
  \dfrac{\dI_k}{\sum_{s=1}^N
    \dI_s-1} \right]
\nonumber\\[0.0ex]
& = \sum_{\ell,r=1}^vR_{i\ell} R_{jr}\sum_{k=1}^N
\frak{E}_\text{diag}\left[\left(X^{(k)}_\ell-\tilde X_\ell\right) \left(X^{(k)}_r-\tilde X_r\right)
  \dfrac{\dI_k}{\sum_{s=1}^N
    \dI_s-1} \right]\,,
\end{align}
with $\frak{E}_\text{diag}$ denoting the sample expectation corresponding to
a multinormal measure with diagonal covariance matrix $\Lambda = \diag(\lambda) =
\trans{R}\Sigma R$, conditioned to $M>1$. Having diagonalized the
product measure, we observe that the integrand on the {\it r.h.s.} is
odd for $\ell\ne r$ and even for $\ell=r$ under the joint change of
variables 
$X^{(k)}_\ell \to - X^{(k)}_\ell$ for $k=1,\ldots,N$, similarly to what we did in sect.~2. As a
consequence, it holds
\begin{equation}
\frak{E}\left[(\hat{\frak{S}}_\cB)_{ij}\right] = \sum_{\ell=1}^v R_{i\ell} R_{j\ell}\left\{\sum_{k=1}^N
\frak{E}_\text{diag}\left[\left(X^{(k)}_\ell-\tilde X_\ell\right)^2\dfrac{\dI_k}{\sum_{s=1}^N
    \dI_s-1} \right]\right\}\,,
\label{eq:scovbias}
\end{equation}
whence we infer that the matrix $\frak{E}[\hat{\frak{S}}_\cB]$ is
diagonalized by the same matrix $R$ as $\Sigma$. From
eq.~(\ref{eq:scovbias}) we also conclude that 
\begin{equation}
\frak{bias}[{\hat{\frak{S}}_\cB}] = R\,\diag(w)\,\trans{R}, \qquad w_i \equiv \sum_{k=1}^N
\frak{E}_\text{diag}\left[\left(X^{(k)}_i-\tilde X_i\right)^2\dfrac{\dI_k}{\sum_{s=1}^N
    \dI_s-1} \right] - \mu_i\,,
\end{equation}
$i=1,\ldots,v$. It should be observed that in general $w_i\ne \frak{bias}[\hat\mu_i]$ since the computation of $\hat\mu_i$ requires
the diagonalization of ${\hat{\frak{S}}_\cB}$, which is in general performed by a diagonalizing matrix
$\hat R\ne R$. Nevertheless, if $w$ vanishes then $\frak{bias}[{\hat{\frak{S}}_\cB}]$ vanishes too. Now, 
we observe that $w_i$ splits into three contributions,
\begin{align}
w_{i1} & = \phantom{-2}\sum_{k=1}^N \frak{E}_\text{diag}\left[\left(X^{(k)}_i\right)^2\dfrac{\dI_k}{\sum_{s=1}^N
    \dI_s-1} \right] - \mu_i\,,\\[0.0ex]
w_{i2} & = -2\sum_{k=1}^N
\frak{E}_\text{diag}\left[X^{(k)}_i\tilde X_i\dfrac{\dI_k}{\sum_{s=1}^N
    \dI_s-1} \right]\,,\\[0.0ex]
w_{i3} & = \phantom{-2}\sum_{k=1}^N \frak{E}_\text{diag}\left[\left(\tilde X_i\right)^2\dfrac{\dI_k}{\sum_{s=1}^N
    \dI_s-1} \right]\,,
\end{align}
which can be exactly calculated and expressed in terms of $\mu_i$, $\alpha$ and $N$. For instance, 
\begin{align}
w_{i1} & = N\frak{E}_\text{diag}\left[\frac{1}{M-1}(X_i^{(1)})^2\dI_1\,\biggr|\,M>1\right] - \mu_i \nonumber\\[1.0ex] 
& = N\sum_{m=2}^N\frac{1}{m-1}\frak{E}_\text{diag}\left[(X_i^{(1)})^2\dI_1\,|\,M=m\right]-\mu_i\nonumber\\[1.0ex]
& = N\sum_{m=2}^N\frac{1}{m-1}\alpha\mu_i{N-1\choose m-1}\alpha^{m-1}(1-\alpha)^{N-m} - \mu_i\nonumber\\[1.0ex]
& = \mu_i\sum_{m=2}^N\frac{m}{m-1}{N\choose m}\alpha^m(1-\alpha)^{N-m} -\mu_i\,.
\end{align}
Analogously, we have
\begin{equation}
w_{i2} = -2\mu_i\sum_{m=2}^N\frac{1}{m-1}{N\choose m}\alpha^m(1-\alpha)^{N-m}\,,
\end{equation}
\begin{equation}
w_{i3} = \phantom{-2}\mu_i\sum_{m=2}^N\frac{1}{m-1}{N\choose m}\alpha^m(1-\alpha)^{N-m}\,.
\end{equation}
Hence, it follows 
\begin{equation}
w_i = -\mu_i[1 + \alpha(N -1)](1-\alpha)^{N-1}\,.
\end{equation}
Since $\alpha>0$, we see that $\lim_{N\to\infty}w_i=0$. 
Thus, we conclude that ${\hat{\frak{S}}_\cB}$ is asymptotically unbiased.

\begin{figure}
  \centering
  \includegraphics[width=1.0\textwidth]{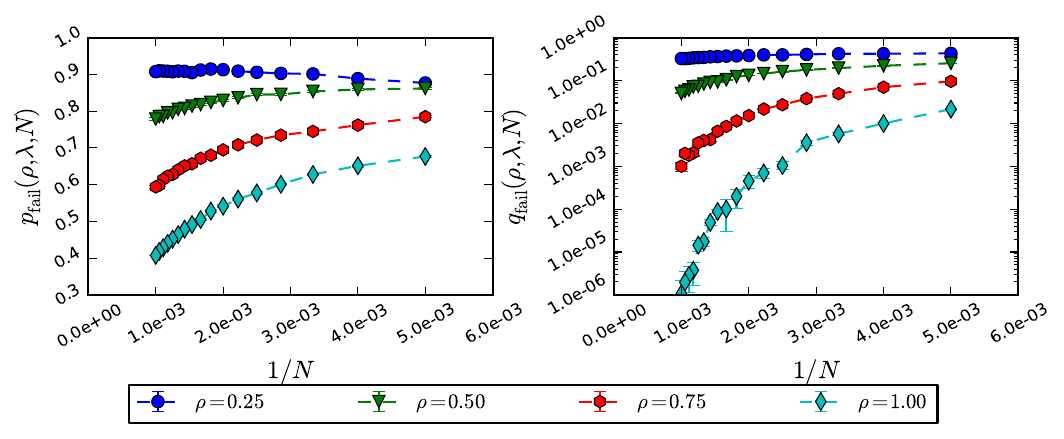}
  \caption{\small {\it Left}: numerical reconstruction of the failure probability of the iterative 
procedure as $v=4$ and $\Sigma=\diag(0.1,0.3,0.8,2.2)$, for several values of $\rho$ and for 
$N=200,250,\ldots,1000$.  {\it Right}: failure probability of the perturbative regularization 
with same parameters.
    \label{fig:pfail}}
\end{figure}

A discussion of the variance of the sample truncated covariance matrix is beyond the scope
of the present paper. We just observe that, apart from the above calculation, studying the sample 
properties of the truncated spectrum is made hard by the fact that eigenvalues and eigenvectors
of a diagonalizable matrix are intimately related from their very definition, thus such study would require
a careful consideration of the distribution of the sample diagonalizing matrix $\hat R \ne R$ of ${\hat{\frak S}}_B$.

\subsection{Perturbative regularization of $\tau_\rho^{-1}$}

When $\mu$ is critically close to the internal boundary of $\cD(\tau_{\rho}^{-1})$, a sample estimate $\hat\mu$ may 
fall outside of it due to statistical fluctuations. In that case the iterative procedure described in the previous sections 
diverges. On the quantitative side, the ill-posedness of the reconstruction problem is measured by the 
failure probability
\begin{equation}
p_{\text{fail}}(\rho,\Sigma,N) = \mathds{P}\left[\,\hat\mu \notin
  \cD(\tau_\rho^{-1}) \ \bigr|
  \ X^{(k)}\sim\cN_v\left(0;\Sigma\right)\,,\ k=1,\ldots,N\,\right]\,,
\end{equation}
which is a highly non--trivial function of $\rho$, $\Sigma$ and $N$. An illustrative example of it is reported 
in Fig.~\ref{fig:pfail} (left), which refers to a specific case with $v=4$ and $\Sigma=\diag(0.1,0.3,0.8,2.2)$. 
The plot suggests that the iterative procedure becomes severely ill--posed in the regime of strong truncation. 

In order to regularize the problem, we propose to go back to eq.~(\ref{eq:truncspectrum}) and consider it from a
different perspective. Specifically, we move  from the observation that a simplified framework occurs in the special
circumstance when the eigenvalue spectra are fully degenerate, which is essentially equivalent to the set--up of 
ref.~\cite{tallis}. If $\mu_1=\ldots = \mu_v \equiv \tilde\mu$, 
by symmetry arguments it follows $\lambda_1 = \ldots = \lambda_v \equiv \tilde\lambda$ and the other way round. 
Eq.~(\ref{eq:truncspectrum}) reduces in this limit to
\begin{equation}
\tilde\mu = \tilde\lambda\frac{F_{v+2}}{F_v}\left(\frac{\rho}{\tilde\lambda}\right) \equiv {\cal T}_\rho(\tilde\lambda)\,,
\label{eq:pertspectrum}
\end{equation}
It can be easily checked that the function ${\cal T}_\rho(\tilde\lambda)$ is monotonic increasing in $\tilde\lambda$. In addition, 
we have
\begin{equation}
i\text{)}\quad \lim_{\tilde\lambda\to 0}{\cal T}_\rho(\tilde\lambda) = 0 \,,
\qquad\qquad ii\text{)}\lim_{\tilde\lambda\to\infty}{\cal T}_\rho(\tilde\lambda) = \frac{\rho}{v+2}\,,
\end{equation}
thus eq.~(\ref{eq:pertspectrum}) can be surely (numerically) inverted provided $0<\tilde\mu<\rho/(v+2)$. 
We can regard eq.~(\ref{eq:pertspectrum}) as an approximation to the original problem, eq.~(\ref{eq:truncspectrum}). 
When $\mu$ is not degenerate, we must define $\tilde\mu$ in terms of the components of $\mu$. 
One possibility is to average them, \ie to choose
\begin{equation}
\tilde\mu = \frac{1}{v}\sum_{i=1}^v\mu_i\, . 
\label{eq:muaver}
\end{equation}
Subject to this, we expect $\tilde\lambda$ to lie somewhere between $\lambda_1$ and $\lambda_v$. 
Eq.~(\ref{eq:pertspectrum}) can be thought of as the lowest order approximation of a perturbative expansion
of eq.~(\ref{eq:truncspectrum}) around the point $\lambda_\text{T}=\{\tilde\lambda,\ldots,\tilde\lambda\}$.
If the condition number of $\Sigma$ is not extremely large, such an expansion is expected to quickly converge, 
so that a few perturbative corrections to $\lambda_\text{T}$ should be sufficient to guarantee a good level of 
approximation. 
\begin{figure}[t!]
  \begin{minipage}[!t]{0.49\textwidth}
    \hskip 0.5cm
    \includegraphics[width=0.8\textwidth]{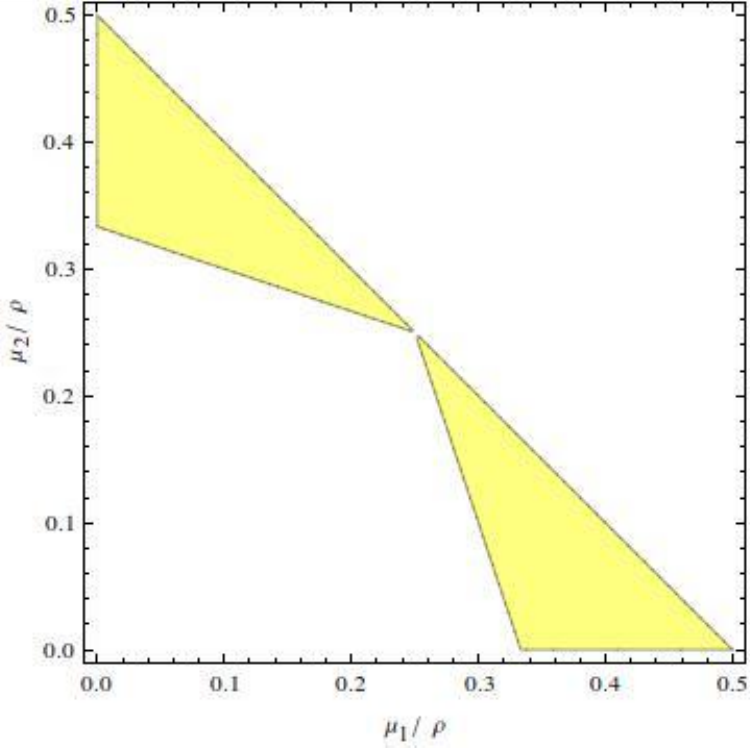}
    \hfill
  \end{minipage}
  \hskip 1.0cm\begin{minipage}[!t]{0.49\textwidth}
    \vskip -1.0cm
    \hfill
    \includegraphics[width=0.8\textwidth]{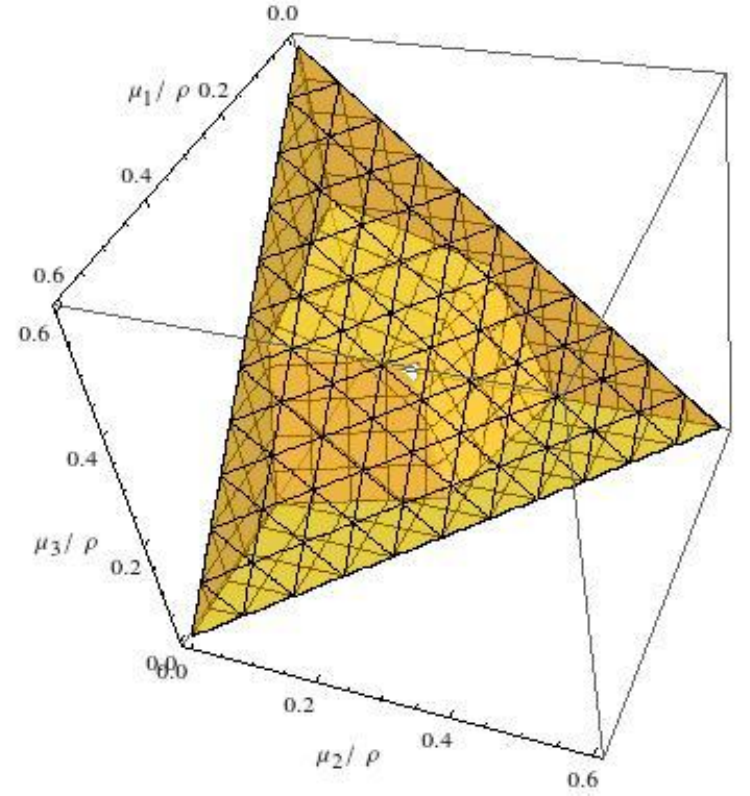}
    \hbox{\hskip 1.0cm \phantom{a}}
  \end{minipage}
  \vskip -0.2cm
  \caption{ \footnotesize {\it Left}: set difference $\cD({\cal T}_\rho^{-1})\setminus \cD(\tau_\rho^{-1})$ in $v=2$ dimensions. {\it Right}:
  set difference $\cD({\cal T}_\rho^{-1})\setminus \cD(\tau_\rho^{-1})$ in $v=3$ dimensions.}
  \label{fig:setdiff}
  \vskip 0.5cm
\end{figure}

As mentioned above, a technical implementation of the perturbative approach and a thorough discussion of its
properties are deferred to a separate paper \cite{palopert}. Here, we limit ourselves to observing that 
the definition domain of perturbation theory is ultimately set by its lowest order approximation, since 
corrections to eq.~(\ref{eq:pertspectrum}) are all algebraically built in terms of it, with no additional
constraints. Following eq.~(\ref{eq:muaver}), the domain of  ${\cal T}_\rho^{-1}$ comes to be defined as
\begin{equation}
\cD({\cal T}_\rho^{-1}) = \left\{\mu\in\dR^v_+:\ \ \sum_{i=1}^v \mu_i \le \frac{\rho v}{v+2}\right\}\,,
\end{equation}
and it is clear that $\cD(\tau_\rho^{-1})\subset\cD({\cal T}_\rho^{-1})$ (it is sufficient to sum term by term all
the inequalities contributing to eq.~(\ref{eq:domain})). In Fig.~\ref{fig:setdiff}, 
we show the set difference $\cD({\cal T}_\rho^{-1})\setminus\cD(\tau_\rho^{-1})$ in $v=2$ and $v=3$ 
dimensions. When $\mu\in\cD(\tau_\rho^{-1})$ but its estimate $\hat\mu\notin\cD(\tau_\rho^{-1})$, it 
may well occur $\hat\mu\in\cD({\cal T}_\rho^{-1})$, \ie the set difference acts as an absorbing shield of
the statistical noise. Therefore, if we define the failure probability of the perturbative 
reconstruction as
\begin{equation}
q_{\text{fail}}(\rho,\Sigma,N) = \mathds{P}\left[\,\hat\mu \notin
  \cD({\cal T}_\rho^{-1}) \ \bigr|
  \ X^{(k)}\sim\cN_v\left(0;\Sigma\right)\,,\ k=1,\ldots,N\,\right]\,,
\end{equation}
we expect the inequality $q_{\text{fail}}(\rho,\Sigma,N)\ll p_{\text{fail}}(\rho,\Sigma,N)$ 
to generously hold. An example is given in Fig.~\ref{fig:pfail} (right): we see that
$q_\text{fail}$ becomes lower than $p_\text{fail}$ by orders of magnitude as soon as $\rho$ and 
$N$ are not exceedingly small. In this sense, the operator ${\cal T}_\rho^{-1}$ can be regarded as
the lowest order approximation of a regularizing operator for $\tau_\rho^{-1}$. 

\section{Conclusions}

In this paper we have studied how to reconstruct the covariance matrix $\Sigma$
of a normal multivariate $X\sim\cN_v(0,\Sigma)$ from the matrix
$\fS_\cB$ of the spherically truncated second moments, describing the covariances 
among the components of $X$ when the probability density is cut off outside 
a centered Euclidean ball. We have shown that $\Sigma$ and $\fS_{\cB}$ share 
the same eigenvectors. Therefore, the problem amounts to relating the eigenvalues 
of $\Sigma$ to those of $\fS_{\cB}$. Such relation entails the inversion of a system 
of non--linear integral equations, which admits unfortunately no closed--form solution. 
Having found a necessary condition for the invertibility of the system, we have shown 
that the eigenvalue reconstruction can be achieved numerically via a converging fixed 
point iteration. In order to prove the convergence, we rely ultimately upon some 
probability inequalities, known in the literature as {\it square correlation 
inequalities}, which have been recently proved in \cite{mukerjee}.

In order to explore the convergence rate of the fixed point iteration, we
have implemented some variations of the non--linear Gauss--Jacobi scheme. Specifically, 
we have found that over--relaxing the basic iteration enhances the convergence rate 
by a moderate factor. However, the over--relaxed algorithm still slows down
exponentially in the number of eigenvalues and polynomially in the truncation
radius of the Euclidean ball. We have shown that a significant reduction of the slowing 
down can be achieved in the regime of strong truncation by adapting the relaxation 
parameter to the eigenvalue it is naturally associated with, so as to boost the 
higher components of the spectrum. 

We have also discussed how the iterative procedure works when the
eigenvalue reconstruction is performed on sample estimates of the
truncated covariance spectrum. Specifically, we have shown that 
the statistical fluctuations make the problem ill--posed. We have sketched a possible 
way--out based on perturbation theory, which is thoroughly discussed in a separate
paper \cite{palopert}. 

A concrete implementation of the proposed approach requires the computation of
a set of multivariate Gaussian integrals over the Euclidean ball. For this, we have
extended to the case of interest a technique, originally proposed by 
Ruben for representing the probability content of quadratic forms of normal 
variables as a series of chi--square distributions.  In the paper, we have shown 
the practical feasibility of the series expansion for the integrals involved in
our computations. 

\section*{Acknowledgements}

We are grateful to A.~Reale for encouraging us throughout all stages of this work, and
to G.~Bianchi for technical support at ISTAT. We also thank
R.~Mukerjee and S.~H.~Ong for promptly informing us about their proof of
eqs.~(\ref{eq:varineq}) and (\ref{eq:covineq}). The computing
resources used for our numerical study and the related technical
support at ENEA have been provided by the CRESCO/ENEAGRID High
Performance Computing infrastructure and its staff \cite{cresco}. CRESCO
({\color{red} C}omputational {\color{red} RES}earch centre on
{\color{red} CO}mplex systems) is funded by
ENEA and by Italian and European research programmes.

\begin{appendices}

\section{Proof of Theorem 5.1}

As already mentioned in sect.~4, the proof follows in the tracks of the original one
of ref.~\cite{ruben4}. We detail the relevant steps for $\alpha_k$, while
for $\alpha_{jk}$ we only explain why it is necessary to distinguish
between equal or different indices and the consequences for either case.

In order to prove eq.~(\ref{eq:ruben1}), we first express $\alpha_k$ in
spherical coordinates, i.e. we perform the change of variable $x=ru$, being $r=||x||$
and $u\in\partial\cB_v(1)$ (recall that $\rd^v x = r^{v-1}\rd r\,\rd u$, with
$\rd u$ embodying the angular part of the spherical Jacobian and the
differentials of $v-1$ angles); then we insert a factor of $1 =
\exp(r^2/2s)\exp(-r^2/2s)$ under the integral sign. Hence, $\alpha_k$ reads  
\begin{equation}
\alpha_k(\rho;\lambda) =
\frac{1}{(2\pi)^{v/2}|\Lambda|^{1/2}}\int_0^{\sqrt{\rho}}\rd r\,
r^{v-1}\frac{r^2}{\lambda_k}\exp\left(-\frac{r^2}{2s}\right)\int_{\partial\cB_v(1)}\rd u\ u_k^2\exp\left(-\frac{Q(u)r^2}{2}\right)\,.
\label{eq:app1}
\end{equation}
The next step consists in expanding the inner exponential in Taylor
series\footnote{in his original proof, Ruben considers a more general set--up,
  with the center of the Euclidean ball shifted by a vector $b\in\RR^v$ from the
  center of the distribution. In that case, the Gaussian exponential looks
  different and must be expanded in series of Hermite polynomials. Here, we work in
  a simplified set--up, where the Hermite expansion reduces to Taylor's.}, {\it viz.}
\begin{equation}
\exp\left(-\frac{Q(u)r^2}{2}\right) \,
=\,\sum_{m=0}^\infty\, \frac{1}{m!}\frac{r^{2m}}{2^m}(-Q)^m\,. 
\end{equation}
This series converges uniformly in $u$. We review
the estimate just for the sake of completeness: 
\begin{equation}
\left|\sum_{m=0}^\infty\, \frac{1}{m!}\frac{r^{2m}}{2^m}(-Q)^m\right| \le
\sum_{m=0}^\infty \frac{1}{m!}\frac{r^{2m}}{2^m}q_0^m = \exp\left(\frac{r^2q_0}{2}\right)\,,
\end{equation}
being $q_0 = \max_{i}|1/s - 1/\lambda_i|$. It follows that we can integrate
the series term by term. With the help of the uniform average operator introduced in 
eq.~(\ref{eq:unifaverop}), $\alpha_k$ is recast to
\begin{equation}
\alpha_k(\rho;\lambda) = \sum_{m=0}^\infty \frac{1}{m!}
\frac{1}{\lambda_k}\frac{1}{|\Lambda|^{1/2}}\M[(-Q)^mu_k^2]\frac{1}{2^{v/2+m-1}\Gamma(v/2)}\int_0^{\sqrt{\rho}}\rd
r\, r^{v+2m+1}\exp\left(-\frac{r^2}{2s}\right)\,.
\label{eq:auxeq}
\end{equation}
The presence of an additional factor of $u_k^2$ in the angular average is
harmless, since $|u_k^2|<1$. We finally notice that the radial integral can
be expressed in terms of a cumulative chi--square distribution function on
replacing $r\to \sqrt{rs}$, namely 
\begin{equation}
\int_0^{\sqrt{\rho}}\rd r\, r^{v+2m+1}\exp\left(-\frac{r^2}{2s}\right) \, = \,
2^{v/2+m}s^{v/2+m+1}\Gamma\left(\frac{v}{2}+m+1\right)F_{v+2(m+1)}\left(\frac{\rho}{s}\right)\,.
\label{eq:radialint}
\end{equation}
Inserting eq.~(\ref{eq:radialint}) into eq.~(\ref{eq:auxeq}) results in  Ruben's
representation of $\alpha_k$. This completes the first part of the proof. 

As a  next step, we wish to demonstrate that the function $\psi_k$ of
eq.~(\ref{eq:genfunck}) is the generating function of the coefficients
$c_{k;m}$. To this aim, we first recall the identities 
\begin{align}
\label{eq:a12}
a^{-1/2} & \, = \, (2\pi)^{-1/2}\int_{-\infty}^{\infty}\rd x\,\phantom{x^2}\exp\left(-\frac{a}{2}x^2\right)\,,\\[1.0ex]
\label{eq:a32}
a^{-3/2} & \, = \, (2\pi)^{-1/2}\int_{-\infty}^{\infty}\rd x\,x^2 \exp\left(-\frac{a}{2}x^2\right)\,,
\end{align}
valid for $a>0$. On setting $a_i=[1-(1-s/\lambda_i)z]$, we see that $\psi_k$
can be represented in the integral form
\begin{align}
\label{eq:A.8}
\psi_k(z) & = \, \left(\frac{s}{\lambda_k}\right)^{3/2}(2\pi)^{-v/2}\int_{-\infty}^{\infty}\rd
x_k\,x_k^2\exp\left(-\frac{1}{2}\left[1-\left(1-\frac{s}{\lambda_k}\right)z\right]x_k^2\right)\nonumber\\[2.0ex]
& \times \, \prod_{i\ne k}\left(\frac{s}{\lambda_i}\right)^{1/2}\int_{-\infty}^{\infty}\rd
x_i\,\exp\left(-\frac{1}{2}\left[1-\left(1-\frac{s}{\lambda_k}
\right)z\right]x_i^2\right)\nonumber
\end{align}
\begin{align}
\phantom{\psi_k(z)} & =\, \frac{s}{\lambda_k}\frac{s^{v/2}}{(2\pi)^{v/2}|\Lambda|^{1/2}}\int_{\RR^v}\rd^v
x\,x_k^2\,\exp\left(-\frac{1}{2}zsQ(x)-\frac{\trans{x}\cdot x}{2}\right)\,,
\end{align}
provided $|z|<\min_i|1-s/\lambda_i|^{-1}$. As previously done, we introduce spherical
coordinates $x=ru$, and expand $\exp\{-\frac{1}{2}zsQ(x)\}
= \exp\{-\frac{1}{2}zsr^2Q(u)\}$ in Taylor series. By the same argument as
above, the series converges uniformly in $u$ (the factor of $zs$ does not
depend on $u$), thus allowing term by term integration. Accordingly, we have 
\begin{equation}
\psi_k(z) =
\frac{s}{\lambda_k}\frac{s^{v/2}}{(2\pi)^{v/2}|\Lambda|^{1/2}}\sum_{m=0}^{\infty}z^m\,\frac{s^m}{2^mm!}
\int_0^\infty\rd r\,r^{v+2(m+1)-1}\re^{-r^2/2}\int_{\partial\cB_v(1)}\rd u\,[-Q(u)]^mu_k^2\,.
\label{eq:A.9}
\end{equation}
We see that the {\it r.h.s.} of eq.~(\ref{eq:A.9}) looks similar to
eq.~(\ref{eq:auxeq}), the only relevant differences being the presence of the
factor of $z^m$ under the sum sign and the upper limit of the 
radial integral. With some algebra, we arrive at 
\begin{equation}
\psi_k(z) = \sum_{m=0}^\infty z^m\left\{\frac{2}{m!}\,\frac{s}{\lambda_k}\,\frac{s^{v/2+m}}{|\Lambda|^{1/2}}\,\frac{\Gamma(v/2+m+1)}{\Gamma(v/2)}\,\M[(-Q)^mu_k^2]\right\}\,.
\end{equation}
The series coefficients are recognized to be precisely those of eq.~(\ref{eq:rcoefs1}).

In the last part of the proof, we derive the recursion fulfilled by the
coefficients $c_{k;m}$. To this aim, the $m^{\rm th}$ derivative of $\psi_k$
has to be evaluated at $z=0$ and then identified with $m!\,c_{k;m}$. The key
observation is that differentiating $\psi_k$ reproduces $\psi_k$ itself, that
is to say 
\begin{equation}
\psi'_k(z) = \Psi_{k}(z)\psi_k(z)\,,
\label{eq:psikfirstder}
\end{equation}
with
\begin{equation}
\Psi_k(z) = \frac{1}{2}\sum_{i=1}^v e_{k;i}\left(1-\frac{s}{\lambda_i}\right)\left[1-\left(1-\frac{s}{\lambda_i}\right)z\right]^{-1}\,,
\end{equation}
and the auxiliary coefficient $e_{k;i}$ being defined as in
eq.~(\ref{eq:auxcoefs1}). Eq.~(\ref{eq:psikfirstder}) lies at the origin of
the recursion. Indeed, from  eq.~(\ref{eq:psikfirstder}) it follows that
that $\psi_k''$ is a function of $\psi_k'$ and $\psi_k$, {\it viz.} $\psi_k''
= \Psi_k'\psi_k + \Psi_k\psi_k'$. Proceeding analogously yields the general formula
\begin{equation}
\psi_k^{(m)}(z) = \sum_{r=0}^{m-1}\binom{m-1}{r}\Psi_k^{(m-r-1)}(z)\,\psi_k^{(r)}(z)\,.
\end{equation}
At $z=0$, this reads 
\begin{equation}
m!\,c_{k;m} = \sum_{r=0}^{m-1}\frac{(m-1)!}{(m-r-1)!\,r!}\Psi_k^{(m-r-1)}(0)\ r!\,c_{k;r}\,.
\end{equation}
The last step consists in proving that
\begin{equation}
\Psi_k^{(m)}(0) = \frac{1}{2}m!g_{k;m+1}\,,
\end{equation}
with $g_{k;m}$ defined as in eq.~(\ref{eq:recurs1}). This can be done
precisely as explained in ref.~\cite{ruben4}.

Having reiterated Ruben's proof explicitly in a specific case, it is now easy
to see how the theorem is extended to any other Gaussian integral.  First of
all, from eq.~(\ref{eq:app1}) we infer that each additional subscript in
$\alpha_{k\ell m\ldots}$ enhances the power of the radial coordinate under the
integral sign by 2 units. This entails a shift in the number of degrees of
freedom of the chi--square distributions in Ruben's expansion, amounting to
twice the number of subscripts. For instance, since $\alpha_{jk}$ has two
subscripts, its Ruben's expansion starts by $F_{v+4}$, independently of
whether $j=k$ or $j\ne k$. In second place, we observe that in order to
correctly identify the generating functions of Ruben's coefficients for a
higher--order integral $\alpha_{k\ell m\ldots}$, we need to take into
account the multiplicities of the indices $k$, $\ell$, $m$,\ldots. As an
example, consider the case of $\psi_{jk}$ ($j\ne k$) and  $\psi_{kk}$. 
By going once more through the argument presented in eq.~(\ref{eq:A.8}), we
see that eqs.~(\ref{eq:a12})--(\ref{eq:a32}) are sufficient to show that
eq.~(\ref{eq:genfuncjk}) is the generating function of $\alpha_{jk}$. By
contrast, in order to repeat the proof for the case of $\psi_{kk}$, we need 
an additional integral identity, namely 
\begin{equation}
a^{-5/2}\, = \, \frac{1}{3}\,(2\pi)^{-1/2}\int_{-\infty}^{+\infty}\rd x\, x^4\exp\left(-\frac{a}{2}x^2\right)\,,
\end{equation}
valid once more for $a>0$. Hence, we infer that $\psi_{kk}$ must depend upon
$\lambda_k$ via a factor of $[1-(1-s/\lambda_k)z]^{-5/2}$, whereas $\psi_{jk}$
($j\ne k$) must depend on $\lambda_j$ and $\lambda_k$ via  factors  of
resp. $[1-(1-s/\lambda_j)z]^{-3/2}$ and $[1-(1-s/\lambda_k)z]^{-3/2}$. The
different  exponents are ultimately responsible for the specific values taken
by the auxiliary coefficients $e_{kk;i}$ and $e_{jk;i}$ of eq.~(\ref{eq:auxcoefs2}).

To conclude, we observe that the estimates of the residuals $\cR_{k;m}$ and $\cR_{jk;m}$, 
presented in sect.~4 without an explicit proof, do not require any further technical 
insight than already provided by ref.~\cite{ruben4} plus our considerations. We leave
them to the reader, since they can be obtained once more in the tracks of the 
original derivation of $\cR_{m}$.

\end{appendices}

\bibliographystyle{hunsrt}
\bibliography{main}

\end{document}